\documentclass[12pt]{smfart}

\usepackage[applemac]{inputenc}
\usepackage[francais]{babel}
\usepackage{smfenum}
\usepackage{smfthm}
\usepackage{amssymb}
\usepackage{eucal}

\newcommand{\Qp}{\mathbf{Q}_p}
\newcommand{\Zp}{\mathbf{Z}_p}
\newcommand{\Fp}{\mathbf{F}_p}
\newcommand{\Cp}{\mathbf{C}_p}
\newcommand{\ZZ}{\mathbf{Z}}

\newcommand{\QQ}{\mathbf{Q}}
\newcommand{\OO}{\mathcal{O}}

\newcommand{\Qpnrhat}{\widehat{\mathbf{Q}}_p^{\mathrm{nr}}}
\newcommand{\Qph}{\mathbf{Q}_{p^h}}

\newcommand{\val}{\operatorname{val}}

\renewcommand{\phi}{\varphi}
\renewcommand{\geq}{\geqslant}
\renewcommand{\leq}{\leqslant} 

\renewcommand{\tilde}{\widetilde}

\newcommand{\vl}{\val_\Lambda}
\newcommand{\vp}{\val_p}

\newcommand{\ve}{\val_{\mathrm{E}}}

\newcommand{\eps}{\varepsilon}

\newcommand{\bst}{\mathbf{B}_{\mathrm{st}}} 
\newcommand{\bcris}{\mathbf{B}_{\mathrm{cris}}} 
\newcommand{\be}{\mathbf{B}_{\mathrm{e}}} 
\newcommand{\bmax}{\mathbf{B}_{\mathrm{max}}} 
 
\newcommand{\btp}{\widetilde{\mathbf{B}}^+_{\mathrm{rig}}} 
\newcommand{\btrig}[1]
{\widetilde{\mathbf{B}}^{\dagger #1}_{\mathrm{rig}}} 
\newcommand{\brig}[2]{\mathbf{B}^{\dagger #1}_{\mathrm{rig} #2}}
\newcommand{\brigplus}[1]{\mathbf{B}^+_{\mathrm{rig} #1}}
\newcommand{\bdr}{\mathbf{B}_{\mathrm{dR}}}  
\newcommand{\bfont}{\mathbf{B}}

\newcommand{\bdag}[1]{\mathbf{B}^{\dagger #1}}
\newcommand{\adag}[1]{\mathbf{A}^{\dagger #1}}

\newcommand{\btdag}[1]{\widetilde{\mathbf{B}}^{\dagger #1}}

\newcommand{\atfont}{\widetilde{\mathbf{A}}}
\newcommand{\atplus}{\widetilde{\mathbf{A}}^+}
\newcommand{\btplus}{\widetilde{\mathbf{B}}^+}

\newcommand{\bplus}{\mathbf{B}^+}
\newcommand{\bt}{\widetilde{\mathbf{B}}}
\newcommand{\et}{\widetilde{\mathbf{E}}}
\newcommand{\etplus}{\widetilde{\mathbf{E}}^+}

\newcommand{\dcris}{\mathrm{D}_{\mathrm{cris}}}
\newcommand{\dst}{\mathrm{D}_{\mathrm{st}}}
\newcommand{\ddr}{\mathrm{D}_{\mathrm{dR}}}

\newcommand{\dsen}{\mathrm{D}_{\mathrm{Sen}}}
\newcommand{\dfont}{\mathrm{D}}
\newcommand{\mfont}{\mathrm{M}}
\newcommand{\dtilde}{\widetilde{\mathrm{D}}}
\renewcommand{\ddag}[1]{\mathrm{D}^{\dagger #1}}

\newcommand{\sat}{\mathrm{Sat}}
\newcommand{\mat}{\mathrm{Mat}}
\newcommand{\fil}{\mathrm{Fil}}
\newcommand{\GL}{\mathrm{GL}}
\newcommand{\calM}{\mathcal{M}}
\newcommand{\NP}{\mathrm{NP}}
\newcommand{\im}{\mathrm{im}}
\newcommand{\rg}{\mathrm{rg}}
\newcommand{\rep}{\mathrm{Rep}}

\author{Laurent Berger}
\address{UMPA ENS Lyon \\
46 allée d'Italie \\
69364 Lyon Cedex 07 \\
France}
\email{laurent.berger@umpa.ens-lyon.fr}
\urladdr{www.umpa.ens-lyon.fr/\~{}lberger/}

\date{Avril 2007}

\title{B-paires et $(\varphi,\Gamma)$-modules}

\alttitle{B-pairs and $(\varphi,\Gamma)$-modules}



\NumberTheoremsIn{subsection}

\begin{document}

\begin{abstract}
On étudie 
la catégorie des $B$-paires $(W_e,W_{dR}^+)$ où
$W_e$ est un $\mathbf{B}_{\mathrm{cris}}^{\varphi=1}$-module 
libre muni d'une action
semi-linéaire et continue de $G_K$ et où $W_{dR}^+$
est un $\mathbf{B}_{\mathrm{dR}}^+$-réseau stable par $G_K$ de 
$\mathbf{B}_{\mathrm{dR}} 
\otimes_{\mathbf{B}_{\mathrm{cris}}^{\varphi=1}} W_e$. Cette catégorie
contient celle des représentations $p$-adiques, et est 
naturellement équivalente à la catégorie de tous les
$(\varphi,\Gamma)$-modules sur l'anneau de Robba.
\end{abstract}

\begin{altabstract}
We study 
the category of $B$-pairs $(W_e,W_{dR}^+)$ where
$W_e$ is a free $\mathbf{B}_{\mathrm{cris}}^{\varphi=1}$-module with a 
semilinear and continuous action of $G_K$ 
and where $W_{dR}^+$
is a $G_K$-stable $\mathbf{B}_{\mathrm{dR}}^+$-lattice in  
$\mathbf{B}_{\mathrm{dR}} 
\otimes_{\mathbf{B}_{\mathrm{cris}}^{\varphi=1}} W_e$. 
This category contains the category of $p$-adic representations
and is naturally equivalent to the category of all 
$(\varphi,\Gamma)$-modules over the Robba ring.
\end{altabstract} 

\maketitle

\setcounter{tocdepth}{2}
\tableofcontents

\setlength{\baselineskip}{18pt}

\section*{Introduction}\label{secintro}

Dans tout cet article, $p$ est un nombre premier, $k$ 
est un corps parfait de caractéristique $p$, et $K$ est une
extension finie totalement ramifiée 
de $K_0=W(k)[1/p]$. On s'intéresse aux
représentations $p$-adiques 
du groupe de Galois $G_K = \mathrm{Gal}
(\overline{K}/K)$. La théorie de Hodge $p$-adique (cf. \cite{F3,F4})
a pour but de décrire certaines de ces représentations, celles
qui {\og proviennent de la géométrie \fg}, en termes d'objets
plus maniables, les $(\phi,N)$-modules filtrés. Le résultat le plus
satisfaisant dans cette direction est le théorème de Colmez-Fontaine,
qui dit que le foncteur $V \mapsto \dst(V)$ réalise une équivalence
de catégories entre la catégorie des représentations $p$-adiques
semi-stables et la catégorie des $(\phi,N)$-modules filtrés admissibles.

Si $D$ est un $\phi$-module filtré qui provient de la cohomologie
d'un schéma propre $X$ sur $\OO_K$, alors le $\phi$-module sous-jacent
ne dépend que de la fibre spéciale de $X$ (c'en est la cohomologie
cristalline) alors que la filtration ne dépend 
que de la fibre générique
(c'est la filtration de Hodge de la
cohomologie de de Rham, dans laquelle se plonge la 
cohomologie cristalline). Si $V=V_{\mathrm{cris}}(D)$ 
et $\be=\bcris^{\phi=1}$, alors on voit
que $\be \otimes_{\Qp} V = (\bcris \otimes_{K_0} D)^{\phi=1}$
ne dépend que de la structure de $\phi$-module de $D$ et de plus,
les $\phi$-modules $D_1$ et $D_2$ sont isomorphes si et seulement
si les $\be$-représentations 
$\be \otimes_{\Qp} V_1$ et $\be \otimes_{\Qp} V_2$ le sont (cf.\
la proposition 8.2 de \cite{F03}). Par ailleurs, 
$\bdr^+ \otimes_{\Qp} V = \mathrm{Fil}^0 (\bdr \otimes_{K_0} D)$ et
les modules filtrés
$K \otimes_{K_0} D_1$ et $K \otimes_{K_0} D_2$ 
sont isomorphes si et seulement
si les $\bdr^+$-représentations  
$\bdr^+ \otimes_{\Qp} V_1$ 
et $\bdr^+ \otimes_{\Qp} V_2$ le sont. 

L'idée de cet article est d'isoler les phénomènes liés
à la {\og fibre spéciale \fg} et à la {\og fibre générique \fg} 
en considérant non pas des représentations $p$-adiques $V$,
mais des \emph{$B$-paires} $W=(W_e,W_{dR}^+)$ où
$W_e$ est un $\be$-module libre muni d'une action
semi-linéaire et continue de $G_K$ et où $W_{dR}^+$
est un $\bdr^+$-réseau stable par $G_K$ de $W_{dR}
= \bdr \otimes_{\be} W_e$. 
Rappelons que les anneaux $\bcris$
et $\bdr$ sont reliés, en plus de l'inclusion 
$\bcris \subset \bdr$, par
la \emph{suite exacte fondamentale} :
$0 \to \Qp \to \bcris^{\phi=1} \to \bdr/\bdr^+ \to 0$.
 
Si $V$ est une
représentation $p$-adique, alors on lui associe la $B$-paire
$W(V)=(\be \otimes_{\Qp} V, \bdr^+ \otimes_{\Qp} V)$
et ce foncteur est pleinement fidèle car $\be \cap \bdr^+ = \Qp$
ce qui fait que $V = W_e(V) \cap W_{dR}^+(V)$.

L'un des principaux outils dont on dispose pour étudier les
représentations $p$-adiques est la théorie des 
$(\phi,\Gamma)$-modules. De fait, on a une équivalence de
catégories entre la catégorie des représentations $p$-adiques
et la catégorie des $(\phi,\Gamma)$-modules étales sur l'anneau
de Robba (en combinant des résultats de Fontaine, 
Cherbonnier-Colmez et Kedlaya). Le premier résultat de cet
article est que l'on peut associer 
à toute $B$-paire $W$
un $(\phi,\Gamma)$-module
$\dfont(W)$ sur l'anneau de Robba, et
que ce foncteur est alors une équivalence de catégories.

\begin{enonce*}{Théorème A}
Le foncteur $W \mapsto \dfont(W)$ réalise une équivalence 
de catégories entre la catégorie
des $B$-paires et la catégorie des 
$(\phi,\Gamma)$-modules sur l'anneau de Robba.
\end{enonce*}

La sous-catégorie pleine des $B$-paires de la forme $W(V)$
correspond à la sous-catégorie pleine des 
$(\phi,\Gamma)$-modules étales. On peut alors se demander à
quoi correspondent les $(\phi,\Gamma)$-modules isoclines.
Si $h \geq 1$ et $a \in \ZZ$ sont premiers entre eux, 
soit $\rep(a,h)$ la catégorie dont 
les objets sont les $\Qph$-espaces vectoriels $V_{a,h}$ 
de dimension finie,
munis d'une action 
semi-linéaire de $G_K$ et d'un Frobenius 
lui aussi semi-linéaire $\phi : V_{a,h} \to V_{a,h}$ qui
commute à $G_K$ et qui vérifie $\phi^h=p^a$.
Si $V_{a,h} \in \rep(a,h)$, alors on pose $W_e(V_{a,h}) 
= (\bcris \otimes_{\Qph} V_{a,h})^{\phi=1}$ et 
$W_{dR}^+(V_{a,h}) = \bdr^+ \otimes_{\Qph} V_{a,h}$.

\begin{enonce*}{Théorème B}
Si $V_{a,h} \in \rep(a,h)$, alors $W(V_{a,h}) =(W_e(V_{a,h}),
W_{dR}^+(V_{a,h}))$ est une $B$-paire et le foncteur $V_{a,h}
\mapsto W(V_{a,h})$ définit une équivalence de catégories entre 
$\rep(a,h)$ et la catégorie des $B$-paires $W$ telles que $\dfont(W)$
est isocline de pente $a/h$.
\end{enonce*}

Si la catégorie des $B$-paires est plus riche que la catégorie
des représentations $p$-adiques, la morale de cet article est
qu'elle est aussi plus maniable. Ceci est dû au fait qu'il est
plus facile de travailler avec tous
les $(\phi,\Gamma)$-modules sur l'anneau de Robba que
de se restreindre à ceux qui sont étales. 
La généralisation du théorème de Colmez-Fontaine 
aux $B$-paires devient alors
un simple exercice d'algèbre linéaire.
Si $D$ est un $(\phi,N)$-module filtré, et si
$W_e(D) = (\bst \otimes_{K_0} D)^{\phi=1,N=0}$ 
et $W_{dR}^+(D) = \fil^0(\bdr \otimes_K D_K)$, alors
$W(D)=(W_e(D),W_{dR}^+(D))$ est une $B$-paire, 
qui est semi-stable en un sens évident.
Le foncteur $D \mapsto W(D)$ réalise alors
une équivalence de catégories entre la catégorie des $B$-paires
semi-stables et la catégorie des $(\phi,N)$-modules filtrés.

Pour retrouver le théorème de Colmez-Fontaine à partir de cet
énoncé, il faut identifier quelles sont les $B$-paires semi-stables
qui proviennent des $(\phi,N)$-modules filtrés admissibles. On
peut faire cela en passant par les $(\phi,\Gamma)$-modules,
c'est ce qui est fait dans \cite{Ber5}, où l'on démontre 
suffisamment pour impliquer aussi le théorème de monodromie 
$p$-adique pour les $B$-paires : toute $B$-paire de de Rham
est potentiellement semi-stable.

Dans la suite de l'article, nous donnons quelques 
applications du théorème A. Par exemple, on montre que pour 
tout $\phi$-module $\dfont$ sur l'anneau de Robba, il existe 
un $\phi$-module étale $\mfont \subset \dfont[1/t]$ tel
que $\mfont[1/t] = \dfont[1/t]$. Comme $\dfont(W)[1/t]$
{\og correspond \fg} à $W_e$, on déduit de cet énoncé 
le résultat suivant, où $H_K= \mathrm{Gal}
(\overline{K}/K(\mu_{p^\infty}))$.

\begin{enonce*}{Théorème C}
Si $W_e$ est une $\be$-représentation de $G_K$, alors il
existe une représen\-tation
$p$-adique $V$ de $H_K$ telle que la restriction de $W_e$ à $H_K$ 
est isomorphe à $\be \otimes_{\Qp} V$.
\end{enonce*}

Enfin, nous donnons (quand $K \subset K_0(\mu_{p^\infty})$) une 
construction des $(\phi,\Gamma)$-modules de \emph{hauteur finie}
(ce sont ceux qui ont une base dans laquelle les matrices de $\phi$
et de $\gamma \in \Gamma_K$ n'ont pas de dénominateurs en $X$)
à partir de {\og $\phi$-modules filtrés sur $K_0$ avec action de
$\Gamma_K$ \fg}. Cela permet d'éclaircir la structure des $B$-paires
de hauteur finie, et donc en particulier des représentations de
hauteur finie. Un \emph{$\phi$-module filtré sur $K_0$ avec 
action de $\Gamma_K$} est un 
$\phi$-module $D$ sur $K_0$ muni d'une action de $\Gamma_K$
commutant à $\phi$ et d'une filtration 
stable par $\Gamma_K$ sur
$D_\infty = K_\infty \otimes_{K_0} D$.

\begin{enonce*}{Théorème D}
Si $D$ est un $\phi$-module filtré sur $K_0$ avec 
action de $\Gamma_K$, alors la
$B$-paire $W(D)=((\bmax \otimes_{K_0} D)^{\phi=1}, 
\bdr^+ \otimes_{K_\infty} D_\infty)$ est de hauteur finie.

De plus, toute $B$-paire de hauteur finie s'obtient 
de cette manière et $W(D_1) = W(D_2)$ si et seulement s'il
existe un isomorphisme $K_0[t,t^{-1}] \otimes_{K_0} D_1
= K_0[t,t^{-1}] \otimes_{K_0} D_2$ compatible à $\phi$
et $\Gamma_K$, et compatible
à la filtration quand on étend les scalaires à
$K_\infty(\!(t)\!)$.
\end{enonce*}

Des objets de nature similaire ont été étudiés, dans un cadre un peu
différent, par Kisin dans \cite{Ki06}.

Pour terminer cette introduction, signalons que
la notion de représentation trianguline de \cite{PCtr,PCalt} s'écrit
de manière très agréable en termes de $B$-paires.
On dit qu'une $B$-paire est \emph{trianguline} si elle 
est extension successive de $B$-paires de dimension $1$.
Par le théorème A, cela revient à dire que $\dfont(W)$ est
extension successive de $(\phi,\Gamma)$-modules de rang $1$. 
On retrouve
donc la définition de \cite{PCtr,PCalt} dans le cas où $W = W(V)$. 

\vspace{\baselineskip}
\noindent\textbf{Remerciements:} 
Je remercie Pierre Colmez pour de nombreuses
discussions sur plusieurs thèmes qui se sont concrétisés
dans ce texte, Jean-Marc
Fontaine pour avoir suggéré que les $\be$-représentations
sont de bons objets et m'avoir expliqué ses idées sur le lien
entre les $B$-paires et les presque $\Cp$-représentations,
Kiran Kedlaya pour ses explications sur les $\phi$-modules
et les nombreux anneaux qui interviennent dans ses constructions, 
et Kiran Kedlaya et Ruochuan Liu pour leur construction d'un contre-exemple
à une de mes affirmations (cf.\ le (2) de la remarque \ref{bvunik}).
Cet article a été rédigé à l'IHES qui m'a offert d'excellentes conditions
de travail.

\section{Rappels et compléments}\label{secprelim}

Dans ce chapitre, nous donnons des rappels sur les anneaux
de périodes et les $(\phi,\Gamma)$-modules. 

\subsection{L'anneau $\btrig{}$ et 
ses petits camarades}\label{pmbtr}

Nous commençons par faire des rappels très 
succints sur les définitions
(données dans \cite{F3} par exemple)
des divers anneaux que nous utilisons dans cet article. Rappelons
que $\etplus = \varprojlim_{x \mapsto x^p} \OO_{\Cp}$ est
un anneau de caractéristique $p$, complet pour la valuation
$\ve$ définie par $\ve(x)= 
\vp(x^{(0)})$ et qui contient un élément $\eps$ tel que
$\eps^{(n)}$ est une racine primitive $p^n$-ième de l'unité.
On fixe un tel $\eps$ dans tout l'article.
L'anneau $\et = \etplus[1/(\eps-1)]$ est alors un corps qui 
contient comme sous-corps dense la clôture algébrique
de $\Fp(\!(\eps-1)\!)$. On pose $\atplus=W(\etplus)$ et 
$\atfont=W(\et)$ ainsi que $\btplus=\atplus[1/p]$ et 
$\bt=\atfont[1/p]$. 
L'application $\theta : \btplus \to \Cp$ qui à
$\sum_{k \gg -\infty} p^k [x_k]$ associe 
$\sum_{k \gg -\infty} p^k x_k^{(0)}$ est un morphisme
d'anneaux surjectif et $\bdr^+$ est le complété de $\btplus$
pour la topologie $\ker(\theta)$-adique, ce qui en fait un espace
topologique de Fréchet. On pose $X=[\eps]-1 \in \atplus$ et
$t=\log(1+X) \in \bdr^+$ et on définit $\bdr$ par $\bdr=\bdr^+[1/t]$.
Soit $\tilde{p} \in \etplus$ un élément tel que $p^{(0)}=p$.
L'anneau $\bmax^+$ est l'ensemble des séries $\sum_{n \geq 0} b_n 
([\tilde{p}]/p)^n$ où $b_n \in \btplus$ et $b_n \to 0$ quand $n \to \infty$
ce qui en fait un sous-anneau de $\bdr^+$ 
muni en plus d'un Frobenius $\phi$ qui est injectif, mais pas surjectif.
On pose $\btp  = \cap_{n\geq 0} \phi^n(\bmax^+)$ ce qui
en fait un sous-anneau de $\bmax^+$ sur lequel $\phi$ est bijectif.
On pose enfin $\bmax=\bmax^+[1/t]$ et $\bst=\bmax[\log [\tilde{p}]]$
ce qui fait de $\bst$ un sous-anneau de $\bdr$ muni de $\phi$ et
d'un opérateur de monodromie $N$. 
Remarquons que l'on travaille souvent avec $\bcris$ plutôt que
$\bmax$ mais le fait de préférer $\bmax$ ne change rien aux
résultats et est plus agréable pour des raisons techniques.

Rappelons que les anneaux $\bmax$
et $\bdr$ sont reliés, en plus de l'inclusion $\bmax \subset \bdr$, par
la \emph{suite exacte fondamentale} $0 \to \Qp \to \bmax^{\phi=1} 
\to \bdr/\bdr^+ \to 0$.
Ce sont ces anneaux que l'on utilise en théorie de Hodge $p$-adique. 
Le point de départ de la théorie des $(\phi,\Gamma)$-modules 
sur l'anneau de Robba (dont
on parle au paragraphe \ref{subphisl}) est la construction d'anneaux
intermédiaires entre $\btplus$ et $\bt$. Si $r>0$, soit $\btdag{,r}$
l'ensemble des $x = \sum_{k \gg -\infty} p^k [x_k] \in \bt$ tels que
$\ve(x_k)+k \cdot pr/(p-1)$ tend vers $+\infty$ quand $k$ augmente. 
On pose $\btdag{} = \cup_{r > 0} \btdag{,r}$, c'est le corps
des éléments surconvergents, défini dans \cite{CC98}.
L'anneau $\btrig{} =  \cup_{r > 0} \btrig{,r}$ 
défini dans \cite[\S 2.3]{Ber1} est en quelque sorte
la somme de $\btp$ et $\btdag{}$; de fait, on a une suite exacte (d'anneaux et 
d'espaces de Fréchet) :
\[ 0 \to \btplus \to \btp \oplus \btdag{,r} \to \btrig{,r} \to 0. \]

Rappelons que $K_0=W(k)[1/p]$; pour $1 \leq n \leq +\infty$, 
on pose $K_n=K(\mu_{p^n})$ et $H_K = \mathrm{Gal}
(\overline{K}/K_\infty)$ et $\Gamma_K=G_K/H_K$.
Si $R$ est un anneau muni d'une action de $G_K$ (c'est le cas
pour tous ceux que nous considérons), on note $R_K = R^{H_K}$.
L'anneau $\btrig{,r}$ contient l'ensemble des séries 
$f(X) = \sum_{k \in \ZZ} f_k X^k$ avec $f_k \in K_0$ telles que
$f(X)$ converge sur $\{ p^{-1/r} \leq |X| < 1 \}$. Cet anneau est
noté $\brig{,r}{,K_0}$. Si $K$ est une extension finie de $K_0$, il
lui correspond par la théorie du corps de normes 
(cf.\ \cite{FW79} et \cite{WI83})
une extension finie
$\brig{,r}{,K}$ qui s'identifie (si $r$ est assez grand) 
à l'ensemble des séries 
$f(X_K) = \sum_{k \in \ZZ} f_k X_K^k$ avec $f_k \in K'_0$ telles que
$f(X_K)$ converge sur $\{ p^{-1/e r} \leq |X_K| < 1 \}$ où 
$X_K$ est un certain élément de $\btdag{}_K$ et $K'_0$ est
la plus grande sous-extension non ramifiée de $K_0$ dans 
$K_\infty$ et $e=[K_\infty : K_0(\mu_{p^\infty})]$.
On pose $\brig{}{,K} = \cup_{r > 0} \brig{,r}{,K}$ et
$\bdag{,r}_K = \brig{,r}{,K} \cap \btdag{}$ et  $\bdag{}_K =
\cup_{r > 0}  \bdag{,r}_K$.
Les anneaux $\btrig{}$ et $\brig{}{,K}$ 
coïncident avec 
les anneaux $\Gamma^{\mathrm{alg}}_{\mathrm{an,con}}$
et $\Gamma_{\mathrm{an,con}}$ définis
dans \cite[\S 2.2]{KK05} (cf.\ en particulier la convention 2.2.16 et la
remarque 2.4.13 de \cite{KK05}). 
Kedlaya les a étudiés en détail et nous rappelons à présent
quelques uns des résultats que nous utilisons dans la suite.

\begin{prop}\label{btpbez}
Les anneaux $\btrig{,r}$ et l'anneau $\btrig{}$ sont de Bézout, ainsi
que les anneaux $\brig{,r}{,K}$ et l'anneau $\brig{}{,K}$. 
Si $R$ est l'un de ces anneaux, et si
$M$ est un sous-$R$-module d'un $R$-module
libre de rang fini, alors les affirmations suivantes sont équivalentes :
\begin{enumerate}
\item $M$ est libre;
\item $M$ est fermé;
\item $M$ est de type fini.
\end{enumerate}
\end{prop}

\begin{proof}
Pour $\btrig{,r}$ ou $\brig{,r}{,K}$, 
c'est le théorème 2.9.6 de \cite{KK05} et pour 
$\btrig{}$ ou $\brig{}{,K}$, 
c'est une conséquence immédiate de ce que $\btrig{}
= \cup_{r > 0} \btrig{,r}$ et $\brig{}{,K} = \cup_{r > 0} 
\brig{,r}{,K}$. 
L'affirmation quant aux sous-modules des modules libres 
est contenue dans le corollaire 2.8.5 
et la définition 2.9.5 de \cite{KK05}.
\end{proof}

\begin{coro}\label{bezmap}
Si $R$ est l'un des anneaux ci-dessus, et si $f : D \to E$ est un morphisme
de $R$-modules libres de rang fini, alors $\im(f)$ et
$\ker(f)$ sont libres de rang fini. De plus, $\ker(f)$ est saturé dans $D$.
\end{coro}

\begin{rema}\label{nobez}
Les anneaux $\btplus$ et $\btp$ ne sont pas de Bézout.
\end{rema}

\begin{proof}
Commençons par $\btplus$.
Soit $\beta_1 > \beta_2  > \cdots$ une suite convexe 
décroissante d'éléments de $\ZZ[1/p]$ 
qui converge vers $\beta >0$; 
si $r \in \ZZ[1/p]$, écrivons $Y^r$ pour $[\tilde{p}^r]$ (rappelons
que $\etplus$ est un anneau parfait et donc que $\tilde{p}^r$ 
est bien déterminé).
On pose $f=\sum_{i \geq 0} p^i Y^{\beta_{2i}}$ et 
$g=\sum_{i \geq 0} p^i Y^{\beta_{2i+1}}$. Supposons 
que l'idéal de $\btplus$ engendré par $f$ et $g$ est principal, 
engendré par un élément $h$. Cet élément est nécessairement
dans $Y^{\geq \beta} \btplus$, puisque $f$ et $g$ le sont et
comme $\atplus/p$ est intègre, on peut supposer que $h \in \atplus$
et que $I = (f,g) \cap \atplus = h \atplus$. On a alors 
$\ve(\overline{h}) \geq \beta$.
Par ailleurs, on a $(f- g Y^{\beta_0-\beta_1})/p \in (f,g)$ 
et cet élément s'écrit $\sum_{i \geq 0} p^i Y^{\beta_{2i+2}} 
(1+O(Y))$ (c'est là qu'on utilise la convexité de la suite)
et en itérant ce procédé, on voit que 
$I$ contient un élément de la forme
$\sum_{i \geq 0} p^i Y^{\beta_{2i+j}} 
(1+O(Y))$ pour tout $j \geq 0$ et donc que
l'image de $I$ dans
$\et^+$ contient des éléments de 
valuation $\beta_j$ pour tout $j \geq 0$.
Ceci entraîne que $\ve(\overline{h}) = \beta$. Si $\beta \not\in \QQ$,
c'est impossible et donc $\btplus$ n'est pas de Bézout.

Si $\btp$ était de Bézout, 
alors il existerait $h \in \btp$ tel que 
$(f,g)= h \btp$ et en utilisant la théorie
des polygones de Newton de \cite[\S 2.5]{KK05}, on montre
que $h \in \btplus$. Comme ci-dessus, on peut supposer que
$h \in \atplus \setminus p \atplus$. Chacun des éléments
$f_j = \sum_{i \geq 0} p^i Y^{\beta_{2i+j}} 
(1+O(Y))$ construits ci-dessus pour $j \geq 0$ 
peut donc s'écrire $f_j = h x_j$ et encore une fois, on a
forcément $x_j \in \atplus$. On a alors 
$\ve(\overline{h}) \leq \ve(\overline{f}_j)
= \beta_j$ ce qui fait que
$\ve(\overline{h}) \leq \beta$.  
Par ailleurs, si $\alpha \in \ZZ[1/p]$ et 
$\alpha < \beta$, alors $Y^\alpha$ divise
$f$ et $g$ ce qui fait que $(f,g) \subset Y^\alpha \btp$
et donc aussi que $h \in Y^\alpha \btp \cap \atplus = 
Y^\alpha \atplus$ 
pour tout $\alpha < \beta$. Si $\beta \not\in \QQ$, alors
$\ve(\overline{h}) < \beta$ et en choisissant 
$\ve(\overline{h}) < \alpha < \beta$, 
on trouve une contradiction,
ce qui fait que $\btp$ n'est pas de Bézout.
\end{proof}

Une grande partie de l'article \cite{KK05} est consacré à l'étude 
des $\phi$-modules sur l'anneau $\btrig{}$ ou sur l'anneau $\brig{}{,K}$.
Nous rappelons à présent quelques résultats quant au premier cas (pour
le deuxième, voir le paragraphe \ref{subphisl}).

\begin{defi}\label{defstd}
Si $h \geq 1$ et $a \in \ZZ$ sont deux entiers premiers entre eux, alors
le \emph{$\phi$-module élémentaire} $M_{a,h}$ est le $\phi$-module 
de base $e_0,\dots,e_{h-1}$ avec $\phi(e_0) = e_1, \dots, \phi(e_{h-2})
= e_{h-1}$ et $\phi(e_{h-1}) = p^a e_0$ 
(cf.\ la définition 4.1.1 de \cite{KK05};
on utilise $(a,h)$ plutôt que $(c,d)$ pour être 
compatible avec les notations
de Colmez et de Fontaine). 
\end{defi}

\begin{prop}\label{dmked}
Si $M$ est un $\phi$-module sur $\btrig{}$, alors il existe des entiers
$a_i,h_i$ tels que $M = \oplus M_{a_i,h_i}$.
\end{prop}

\begin{proof}
Etant donnée la définition 4.5.1 de \cite{KK05}, 
c'est le (a) du théorème 4.5.7 de \cite{KK05}.
\end{proof}

Remarquons que la décomposition n'est pas canonique, mais que l'ensemble des pentes $s_i = a_i / h_i$ comptées avec multiplicités est canonique (cf.\ le (c) du théorème 4.5.7 de \cite{KK05}). Les rationnels que l'on obtient ainsi  sont les \emph{pentes} de $M$.

\begin{coro}\label{umpsurj}
Si $M$ est un $\phi$-module sur $\btrig{}$, alors $1-\phi : M[1/t] \to 
M[1/t]$ est surjective.
\end{coro}

\begin{proof}
Si $n \geq 0$, alors $(1-\phi)(t^{-n} x) = t^{-n} (1-p^{-n} \phi)(x)$
et il suffit donc de montrer que si $n \gg 0$, alors $1-p^{-n} \phi : 
M \to M$ est surjectif. Etant donnée la proposition \ref{dmked} ci-dessus
et le fait que $M_{a,h}(-n) = M_{a-nh,h}$, 
c'est une conséquence immédiate du (b) de la 
proposition 4.1.3 de \cite{KK05}. 
\end{proof}

L'anneau $\bmax^{\phi=1}$ (qui est égal à $\bcris^{\phi=1}$) 
occupe une place centrale dans cet 
article; il est traditionnellement noté $\be$. Etant donnée la
définition de $\btp$, il est clair que l'on a $\be = (\btp[1/t])^{\phi=1}$. 

\begin{lemm}\label{phiregb}
On a $\be = (\btrig{}[1/t])^{\phi=1}$.
\end{lemm}

\begin{proof}
Si $x \in \be$, alors $x \in (\btp[1/t])^{\phi=1} 
\subset (\btrig{}[1/t])^{\phi=1}$. 
Réciproquement, si $x \in (\btrig{}[1/t])^{\phi=1}$, alors
il existe $n \geq 0$ tel que $x = t^n x_n$ avec $x_n \in 
(\btrig{})^{\phi=p^n}$ et le lemme suit alors de la proposition 3.2 de
\cite{Ber1} qui nous dit que $(\btrig{})^{\phi=p^n} =
 (\btp)^{\phi=p^n}$.
\end{proof}

\begin{lemm}\label{unitbe}
On a $\be^\times = \Qp^\times$ et si $z\in\be$ engendre un $\be$-module de rang $1$ stable par $G_K$, alors $z \in \Qp^\times$.
\end{lemm}

\begin{proof}
Soient $x$ et $y$ dans $\be$ tels que $xy=1$ et $v(x)$ et $v(y)$ leurs 
valuations $t$-adiques dans $\bdr$. On a $v(x)+v(y)=0$ et comme
$\be \cap \fil^1 \bdr = 0$, cela entraîne que $v(x)=v(y)=0$. Le fait
que $\be^\times = \Qp^\times$ suit alors du fait que $\be \cap
\fil^0 \bdr = \Qp$.

Si $z\in\be$ engendre un $\be$-module de rang $1$ stable par $G_K$, alors $g(z)/z \in \Qp^\times$ si $g \in G_K$ et l'application $g \mapsto g(z)/z$ 
définit un caractère cristallin de $G_K$. Un résultat classique 
(cf.\ le lemme 5.1.3 de \cite{F4}) dit qu'un tel
$z \in \bmax$ est de la forme $t^n z_0$ avec $n \in \ZZ$ et $z_0 \in 
\Qpnrhat$ et si en plus $\phi(z)=z$, alors $z \in \Qp^\times$.
\end{proof}

En première approximation, on peut d'ailleurs penser à $\be$ comme à l'anneau des polynômes $P(Y) \in \Cp[Y]$ tels que $P(0) \in \Qp$.

\begin{prop}\label{bebez}
L'anneau $\be$ est de Bézout.
\end{prop}

\begin{proof}
Il suffit de montrer que si $f,g \in \be$, alors l'idéal qu'ils engendrent 
est principal. Soit $n \geq 0$ tel que $f, g \in t^{-n} 
(\btp)^{\phi=p^n}$. Comme $\btrig{}$ est de Bézout, il existe 
$h \in \btrig{}$ tel que $t^n f \btrig{} + t^n g \btrig{} = h \btrig{}$.
En particulier, il existe $\alpha, \beta \in \btrig{}$ tels que $t^n f \alpha
+ t^n g \beta =h$. En appliquant $\phi^{\pm 1}$ à cette relation, on trouve que $h$ et $\phi(h)$ engendrent le même idéal de $\btrig{}$. 
Par la proposition 3.3.2 de \cite{KK05}, on peut (quitte à multiplier $h$
par une unité) supposer que $\phi(h) = p^m h$ avec $m \in \ZZ$; 
on en déduit que l'idéal de $\btrig{}[1/t]$ engendré par $f$ et $g$ est engendré par $h/t^m \in \be$. 

Reste à voir que $h/t^m$ est dans l'idéal de $\be$ engendré par $f$ et $g$. On se ramène à montrer que si $f, g \in \be$ engendrent 
$\btrig{}[1/t]$, alors ils engendrent $\be$. Soit $n \geq 0$ tel que
$t^n f, t^n g \in \btrig{}$ et $I$ l'idéal (principal) de $\btrig{}$ qu'ils engendrent. Considérons la suite exacte $0 \to M \to t^n f \btrig{} \oplus t^n g \btrig{} \to 
I \to 0$. Le $\btrig{}$-module $M$ est fermé et donc 
libre de rang fini (égal à $1$) et 
c'est un $\phi$-module. En tensorisant 
la suite exacte par $\btrig{}[1/t]$ 
(qui est plat sur $\btrig{}$) et en 
prenant les invariants par $\phi$, 
on trouve un morceau de suite 
exacte : $f \be \oplus g \be \to 
\be \to M[1/t] / (1 - \phi)$ et 
on conclut par le corollaire \ref{umpsurj}.
\end{proof}

\begin{rema}\label{topobe}
L'anneau $\be$ est la réunion pour $n \geq 0$ des $t^{-n} 
(\bmax^+)^{\phi=p^n}$ et 
chacun d'entre eux est un espace de
Banach. La topologie de $\be$ est celle de la limite
inductive, ce qui en fait un espace LF.
\end{rema}

Si $h \geq 1$, soit $t_h \in \bmax^+$ l'élément construit dans
\cite[\S 9]{CEV} (c'est une période d'un groupe de Lubin-Tate
associé à l'uniformisante $p$ de $\Qph$)
et dont les propriétés sont rappelées au début
de \cite[\S 2.4]{PC03}. On a notamment $\phi^h(t_h) = p t_h$ 
et $\prod_{j=0}^{h-1} \phi^j(t_h) \in \Qp^\times \cdot t$.

\begin{lemm}\label{bphd}
Si $h \geq 1$ et $a \in\ZZ$, alors $(\btp[1/t])^{\phi^h=p^a}$ est
un $\be$-module libre de rang $h$.
\end{lemm}

\begin{proof}
L'application $x \mapsto x/t_h^a$ est un isomorphisme de
$(\btp[1/t])^{\phi^h=p^a}$ sur $(\btp[1/t])^{\phi^h=1}$.
On se ramène donc au cas $a=0$.
Si $\omega\in\Qph$ engendre une base normale
de $\Qph$ sur $\Qp$, 
alors la matrice 
$(\phi^{i+j}(\omega))_{1 \leq i,j \leq h}$ est inversible et l'application
de $(\btp[1/t])^{\phi^h=1}$ dans $\Qph \otimes_{\Qp} \be$ donnée par
$x \mapsto \sum_{i=0}^{h-1} \phi^i(\omega) \otimes \phi^i(x)$ est 
un isomorphisme, d'où le lemme.
\end{proof}

L'anneau $\btrig{,r}$ est muni d'une topologie de Fréchet, 
qui est définie par des valuations $V_{[r;s]}$ avec $s>r$. On
appelle $\bt^{[r;s]}$ le complété de $\btrig{,r}$ pour la valuation
$V_{[r;s]}$ et
$\bfont^{[r;s]}_K$ le complété de $\brig{,r}{,K}$.
Les valuations $V_{[r;s]}$ et les anneaux 
$\bt^{[r;s]}$ ont été définis dans \cite[\S 2.1]{Ber1}
(où ils sont notés $\bt_{[r;s]}$)
et étudiés dans \cite{PC03,KK05} (entre autres) mais il faut faire
attention au fait que les notations sont différentes. Les valuations
sont indexées par des intervalles et si l'on pose $\rho(r) = (p-1)/pr$, 
alors notre intervalle $[r;s]$ coïncide avec l'intervalle 
$[\rho(s);\rho(r)]$ de \cite{PC03} et de \cite{KK05}.

Voici un tableau récapitulatif de quelques unes des notations : 

\renewcommand{\arraystretch}{1.6}
\begin{center}
\begin{tabular}{|ccc|}
\hline
Berger  \cite{Ber1} & Colmez \cite{PC03} & Kedlaya \cite{KK05} \\
\hline
$\btrig{}$ &  $\btrig{}$
&  $\Gamma^{\mathrm{alg}}_{\mathrm{an,con}}$ \\
$\btrig{,\rho(r)}$ &  $\bt^{]0;r]}$
&  $\Gamma^{\mathrm{alg}}_{\mathrm{an,r}}$ \\
$\bt_{[\rho(s);\rho(r)]}$ & --- & $\Gamma^{\mathrm{alg}}_{[r;s]}$ \\
$\btdag{,\rho(r)}$ &  $\bt^{(0;r]}$
&  $\Gamma^{\mathrm{alg}}_r[1/p]$ \\
$\brig{,\rho(r)}{,K}$ & $\bfont^{]0;r]}_K$ & 
$\Gamma_{\mathrm{an},r}$ \\
$\bdag{,\rho(r)}_K$ & $\bfont^{(0;r]}_K$ & 
$\Gamma_r[1/p]$ \\
\hline
\end{tabular} 
\end{center}
\renewcommand{\arraystretch}{1}

L'anneau $\bt^{[r;s]}$ est muni d'une action de $G_K$ et la méthode de
Sen (cf.\ \cite{PC01} et \cite{BC}) permet de simplifier grandement
l'étude des $\bt^{[r;s]}$-représentations de $G_K$.

\begin{prop}\label{tsbrs}
L'anneau $\widetilde{\Lambda} = \bt^{[r;s]}$ vérifie les conditions (TS1), (TS2) et (TS3) avec $\Lambda_{H_K,n}  = 
\phi^{-n}(\bfont^{[p^n r;p^n s]}_K)$ et $\vl = V_{[r;s]}$, 
les constantes $c_1>0$, $c_2>0$ et 
$c_3 > 1/(p-1)$ pouvant être choisies arbitrairement.
\end{prop}

\begin{proof}
Ceci est démontré dans \cite{PC03} : la condition (TS1) résulte 
du lemme 10.1, la condition (TS2) résulte de la proposition 8.12 et
du fait que $\bt^{]0;\rho(r)]}_K$ est dense dans $\bt^{[r;s]}_K$ 
et la condition (TS3) résulte de la proposition 9.10 et de la même densité. 
\end{proof}

\begin{lemm}\label{intri}
Si $r \gg 0$ et si $I$ est un intervalle contenu dans $[r;+\infty[$, alors
$\bfont^I_K \cap \btrig{,r} = \brig{,r}{,K}$.
\end{lemm}

\begin{proof}
Le corollaire 2.5.7 de \cite{KK05} 
nous dit que $x \in  \bfont^J_K$ pour tout intervalle
$I \subset J \subset [r;+\infty[$ ce qui fait que $x \in \brig{,r}{,K}$.
\end{proof}

\subsection{Les $(\phi,\Gamma)$-modules sur 
l'anneau de Robba}\label{subphisl}
Nous donnons quelques rappels et compléments concernant 
les résultats de Kedlaya 
(en essayant de renvoyer systématiquement à 
\cite{KK05} quand c'est possible)
sur les pentes des 
$\phi$-modules sur $\brig{}{,K}$.
Rappelons qu'un $\phi$-module sur un anneau $R$ est un 
$R$-module libre $\dfont$ muni d'un Frobenius $\phi$
tel que $\phi^*(\dfont)=\dfont$. Si cet anneau est en plus muni d'une
action de $\Gamma_K$, alors un $(\phi,\Gamma)$-module
est un $\phi$-module muni d'une action semi-linéaire 
et continue de $\Gamma_K$ qui commute à $\phi$.
Nous renvoyons à l'article de Kedlaya pour la notion de pentes
des $\phi$-modules. Contentons-nous de dire que si $\dfont$
est un $\phi$-module sur $\brig{}{,K}$ alors ses \emph{pentes}
sont les rationnels qui sortent de la proposition \ref{dmked}
appliquée à $\btrig{} \otimes_{\brig{}{,K}} \dfont$.

\begin{lemm}\label{unis}
Si $\dfont$ est un $\phi$-module de rang $1$ sur $\brig{}{,K}$, 
alors la pente de $\dfont$ appartient à $\ZZ$.
\end{lemm}

\begin{proof}
Cette pente appartient à l'image par $\vp$ du corps des coefficients
de $\brig{}{,K}$ et le lemme suit du fait que ce corps est une
extension non ramifiée de $K_0$. On peut aussi dire que la pente
de $M_{a,1}$ est $a\in\ZZ$.
\end{proof}

\begin{defi}\label{defpt}
Si $\dfont$ est un $\phi$-module, son \emph{degré} $\deg(\dfont)$ est
la pente de $\det(\dfont)$.
\end{defi}

\begin{lemm}\label{extetal}
Si $s \in \QQ$, alors un $\phi$-module sur $\brig{}{,K}$
qui est une extension de $\phi$-modules isoclines 
de pente $s$ est lui-même isocline de pente $s$.
\end{lemm}

\begin{proof}
Ecrivons $s=n/d$; si on a une suite exacte $0 \to \dfont_1 \to \dfont \to \dfont_2 \to 0$
avec $\dfont_1$ et $\dfont_2$ isoclines de pente $s$, alors on peut écrire $\dfont_i =
\brig{}{,K} \otimes_{\bdag{}_K} \dfont_i^\dagger$ où $\dfont_i^\dagger$ est 
un $\phi$-module sur $\bdag{}_K$ admettant une base dans laquelle
$p^{-n} \mat(\phi^d) \in \GL_d(\adag{}_K)$. Le lemme résulte alors
de la proposition 7.4.1 de \cite{KK05} qui nous permet de trouver
une base de $\dfont$ dans laquelle
$p^{-n} \mat(\phi^d) \in \GL_d(\adag{}_K)$.
\end{proof}

Rappelons le résultat principal de \cite{KK04} (redémontré et généralisé
dans \cite{KK05} et \cite{KK06}).

\begin{theo}\label{slopefil}
Si $\dfont$ est un $\phi$-module sur $\brig{}{,K}$, alors il existe une unique filtration $0 = \dfont_0 \subset \dfont_1 \subset \cdots \subset \dfont_\ell = \dfont$ par des sous-$\phi$-modules saturés, telle que :
\begin{enumerate}
\item pour tout $1 \leq i \leq \ell$, le quotient $\dfont_i / \dfont_{i-1}$ est isocline;
\item si l'on appelle $s_i$ la pente de $\dfont_i / \dfont_{i-1}$, alors $s_1 < s_2 < \cdots < s_\ell$.
\end{enumerate}
\end{theo}

\begin{rema}\label{slopecalc}
Plaçons-nous dans la situation du théorème \ref{slopefil} 
ci-dessus; par la proposition \ref{dmked}, il existe des entiers
$a_j,h_j$ tels que $\btrig{} \otimes_{\brig{}{,K}} \dfont
= \oplus M_{a_j,h_j}$. Cette décomposition n'est pas canonique,
mais l'ensemble des $a_j/h_j$ est bien défini et coïncide, si l'on
compte les multiplicités, avec l'ensemble des $s_i$ (cf.\ le 
corollaire 6.4.2 de \cite{KK05}).
\end{rema}

\begin{defi}\label{defnp}
Si $\dfont$ est un $\phi$-module de rang $d$ sur $\brig{}{,K}$, 
alors (cf.\ \cite[p.\ 157]{KK04}) 
son \emph{polygone de Newton} $\NP(\dfont)$ est la réunion des segments
d'extrémités $(i,y_i)$ et $(i+1,y_{i+1})$ pour $0 \leq i \leq d-1$
où $y_0=0$ et $y_{i+1}-y_i$ est 
la $i+1$-ième plus petite pente 
de $\dfont$ en comptant les multiplicités. 
Dans les notations du théorème \ref{slopefil}, $\NP(\dfont)$
est la réunion des $\ell$ segments de longueur 
$\rg(\dfont_i/\dfont_{i-1})$ et de
pente $s_i$ pour $1 \leq i \leq \ell$. 
\end{defi}

Remarquons que par le lemme \ref{unis},
les sommets de  $\NP(\dfont)$  sont à coordonnées entières. 
Le dernier sommet
est de coordonnées $(\rg(\dfont),\deg(\dfont))$.

\begin{lemm}\label{above}
Si $\dfont$ est un $\phi$-module sur $\brig{}{,K}$ et si $\mfont$ est un
sous-$\phi$-module de $\dfont$ (pas nécessairement saturé 
ni de même rang), alors 
$\NP(\mfont)$ est au-dessus de $\NP(\dfont)$, 
et si $\NP(\mfont)$ et $\NP(\dfont)$ ont même extrémité, alors
$\mfont=\dfont$.
\end{lemm}

\begin{proof}
Si l'on appelle $\sigma_1,\dots,\sigma_d$ et $\tau_1,\dots,\tau_m$
les pentes de $\dfont$ et $\mfont$ prises avec multiplicité, l'affirmation 
{\og $\NP(\mfont)$ est au-dessus de $\NP(\dfont)$ \fg} est équivalente 
au fait que pour tout $1 \leq k \leq m$, on a $\sigma_1+\sigma_2
+\cdots+\sigma_k \leq \tau_1+\tau_2+\cdots+\tau_k$, ce
qui revient à dire que pour tout $1 
\leq k \leq m$, la plus petite pente
de $\wedge^k \dfont$ est inférieure ou 
égale à la plus petite pente de $\wedge^k \mfont$. 
Ceci suit, après extension des scalaires à $\btrig{}$, du (a) de
la proposition 4.5.14 de \cite{KK05}. 

Enfin si $\NP(\mfont)$ et $\NP(\dfont)$ ont même extrémité, 
alors $\mfont$ et $\dfont$ ont même rang (le 
rang étant la $x$-longueur du polygone de Newton) et $\det(\mfont) = 
\alpha \cdot \det(\dfont)$ où $\alpha \cdot \brig{}{,K}$ est étale,
$\NP(\mfont)$ et $\NP(\dfont)$ ayant la même extrémité.
Par la proposition 3.3.2 de \cite{KK05}, $\alpha$ est une unité de
$\btrig{}$; c'est donc une unité de $\brig{}{,K}$,
ce qui fait que $\mfont=\dfont$. 
\end{proof}

\begin{coro}\label{ptext}
Si $0 \to \dfont_1 \to \dfont \to \dfont_2$ est une suite exacte de $\phi$-modules,
avec $\rg(\dfont)=\rg(\dfont_1)+\rg(\dfont_2)$,
alors $\deg(\dfont) \geq \deg(\dfont_1) + \deg(\dfont_2)$ 
et on a égalité si et seulement
si $\dfont \to \dfont_2$ est surjective.
\end{coro}

\begin{proof}
Par le lemme \ref{above} ci-dessus, on a $\deg(\im(\dfont \to \dfont_2))
\geq \deg(\dfont_2)$ avec égalité si et seulement
si $\dfont \to \dfont_2$ est surjective. On se ramène donc à montrer que
si on a une suite exacte $0 \to \dfont_1 \to \dfont \to \dfont_2 \to 0$, alors 
$\deg(\dfont) = \deg(\dfont_2) + \deg(\dfont_1)$. Ceci suit du fait que 
$\det \mat(\phi \mid \dfont) = \det \mat(\phi \mid \dfont_1) 
\cdot \det \mat(\phi \mid \dfont_2)$.
\end{proof}

Pour terminer ce paragraphe, 
disons quelques mots de la 
{\og localisation 
en $\eps^{(n)}-1$ \fg}; si $n \geq 0$, alors on pose $r_n=p^{n-1}(p-1)$.
Si $x=\sum_{k \gg -\infty} p^k
[x_k] \in \btdag{,r_0}$ alors la série qui définit $x$ converge 
dans $\bdr^+$ et on en déduit un morphisme injectif noté
$\iota_0$ de $\btdag{,r_0}$ dans $\bdr^+$. Ce morphisme 
s'étend par continuité en un morphisme toujours injectif 
(voir la proposition 2.25 de \cite{Ber1}) $\iota_0 : 
\btrig{,r_0} \to \bdr^+$.
Si $r>0$, alors soit $n(r)$ 
le plus petit entier tel que $r_n \geq r$.
Si $n \geq n(r)$, et si $x \in \btrig{,r}$, alors 
$\phi^{-n}(x) \in  \btrig{,r/p^n} \subset \btrig{,r_0}$
et on en déduit un morphisme injectif $\iota_n = \iota_0
\circ \phi^{-n} : \btrig{,r} \to \bdr^+$.
On pose $Q_1=\phi(X)/X$ 
et $Q_n=\phi^{n-1}(Q_1)$ 
ce qui fait que si $n \geq 1$, alors $Q_n$
est le polynôme minimal de $\eps^{(n)}-1$. Par la
proposition 4.8 de \cite{Ber1}, l'application $\theta 
\circ \iota_n : \btrig{,r} \to \Cp$ est surjective
et son noyau est l'idéal engendré par $Q_n$.
Par le lemme 5.11 de \cite{Ber1}, la 
restriction de $\iota_n$ à $\brig{,r}{,K}$ 
a pour image un sous-anneau dense 
(pour la topologie $t$-adique)
de $K_n[\![t]\!]$.
Si $\dfont^r$ est un $\brig{,r}{,K}$-module, 
cette application $\iota_n : \brig{,r}{,K} \to K_n[\![t]\!]$
nous permet de définir la localisation 
$K_n[\![t]\!] \otimes_{\brig{,r}{,K}}^{\iota_n}
\dfont^r$ de $\dfont^r$
en $\eps^{(n)}-1$ pour $n \geq n(r)$.
Cette construction est fondamentale autant dans cet article
que dans \cite{Ber5}.

\section{Les $B$-paires}\label{secberep}

Une \emph{$B$-paire} est un couple 
$W = (W_e,W_{dR}^+)$ où $W_e$ est 
un $\be$-module libre de rang fini muni 
d'une action semi-linéaire et continue de $G_K$ 
(c'est-à-dire une \emph{$\be$-représentation}) et
$W_{dR}^+$ est un $\bdr^+$-réseau de 
$W_{dR} = \bdr \otimes_{\be} W_e$ stable par $G_K$ 
(rappelons que toute $\bdr$-représentation 
admet un $\bdr^+$-réseau stable par Galois, 
cf.\ par exemple le début du \S 3.5 de \cite{F04}). 

L'objet de ce chapitre est d'étudier la catégorie des $B$-paires, et
le résultat principal est que cette catégorie est  
équivalente à la catégorie des 
$(\phi,\Gamma)$-modules sur $\brig{}{,K}$.

\subsection{La catégorie des $B$-paires}\label{subdef}
Ce paragraphe est principalement consacré à donner des définitions
relatives à la catégorie des $B$-paires.

\begin{defi}\label{defmor}
Si $W$ est une $B$-paire, alors on 
appelle \emph{dimension de $W$} 
le rang commun de $W_e$ et 
de $W_{dR}^+$.

Si $W$ et $X$ sont deux $B$-paires, un \emph{morphisme 
de $B$-paires} $f : W \to X$ est la donnée de deux applications $(f_e,f_{dR}^+)$ de $W_e$ dans $X_e$ et de $W_{dR}^+$ dans $X_{dR}^+$ telles que les deux applications induites par extension des scalaires à $\bdr$ coïncident; on appelle alors $f_{dR}$ cette application.
\end{defi}

\begin{exem}\label{exbr}
Voici deux classes importantes de $B$-paires.
\begin{enumerate}
\item Si $V$ est une représentation $p$-adique de $G_K$, alors 
$(\be \otimes_{\Qp} V, \bdr^+ \otimes_{\Qp} V)$ 
est une $B$-paire notée $W(V)$;
\item si $D$ est un $(\phi,N)$-module filtré sur $K$, alors 
$((\bst \otimes_{K_0} D)^{\phi=1,N=0}, 
\fil^0(\bdr \otimes_K D_K))$ est une $B$-paire notée $W(D)$.
\end{enumerate}
\end{exem}

\begin{lemm}\label{dimun}
Si $W$ est une $B$-paire de dimension $1$, alors il existe un
caractère $\eta : G_K \to \Qp^\times$ et $i \in \ZZ$ tels que
$W_e = \be(\eta)$ et $W_{dR}^+ = t^i \bdr^+$. 
\end{lemm}

\begin{proof}
La première assertion résulte du lemme \ref{unitbe} et la deuxième du
fait que les idéaux fractionnaires de $\bdr$ sont
tous de la forme $t^i \bdr^+$.
\end{proof}

\begin{lemm}\label{satim}
Si $f : W \to X$ est un morphisme de $B$-paires, alors 
$X_e / f_e(W_e)$ est sans torsion.
\end{lemm}

\begin{proof}
Le $\be$-module $\sat f_e(W_e)$ est libre de même rang que $f_e(W_e)$ 
et le lemme \ref{unitbe} montre que $\det f_e(W_e) = \det \sat f_e(W_e)$ ce qui fait que l'image de $f_e(W_e)$ est saturée dans $X_e$.
\end{proof}

\begin{rema}\label{strmor}
En revanche, $X_{dR}^+ / f_{dR}^+(W_{dR}^+)$ peut avoir de la
torsion (considérer par exemple l'application naturelle $(\be,t \bdr^+) \to
(\be, \bdr^+)$). On dit que $f$ est \emph{strict} si  
$X_{dR}^+ / f_{dR}^+(W_{dR}^+)$ est sans torsion. L'existence de
morphismes non-stricts implique que la $\otimes$-catégorie additive des
$B$-paires n'est pas abélienne.
\end{rema}

\begin{defi}\label{ssrep}
Si $W$ et $X$ sont deux $B$-paires, on dit que $W$ est un
\emph{sous-objet} de $X$ s'il existe un morphisme injectif strict
$f : W \to X$. Une \emph{suite exacte} est une suite 
où les morphismes sont stricts et les conditions 
habituelles (image - noyau) sont satisfaites.
\end{defi}

\begin{lemm}\label{morbp}
Si $f : W \to X$ est un morphisme de $B$-paires, alors
$\ker(f) = (\ker (f_e), \ker (f_{dR}^+))$ est un sous-objet
de $W$ et $\im(f) = (\im (f_e), \im (f_{dR}^+)) \subset X$ 
est une $B$-paire et la suite $0 \to \ker(f) \to
W \to \im(f) \to 0$ est exacte. 
\end{lemm}

\begin{defi}\label{modifrep}
Si $W$ et $X$ sont deux $B$-paires, on dit que $W$ est une
\emph{modification} de $X$ si $W_e \simeq X_e$.
\end{defi}

\subsection{Construction de $(\phi,\Gamma)$-modules}\label{subpg}
Dans ce paragraphe, on associe à toute $B$-paire un 
$(\phi,\Gamma)$-module et on montre le théorème A
de l'introduction. 
Les constructions qui permettent de relier la catégorie des $B$-paires 
et celle des $(\phi,\Gamma)$-modules sont proches de celles qu'on trouve 
dans l'article \cite{PC03}, notamment les \S\S 2, 3 et 10.

Rappelons que si $n \geq 0$, alors on pose $r_n=p^{n-1}(p-1)$
et que si $r>0$, alors $n(r)$ est le plus petit entier tel que
$r_n \geq r$. Si $W$ est une $B$-paire est si $n \geq n(r)$,
alors l'application $\iota_n : \btrig{,r}[1/t] \to \bdr$ nous
donne un morphisme $\iota_n : 
\btrig{,r}[1/t] \otimes_{\be} W_e 
\to \bdr \otimes_{\be} W_e = W_{dR}$.

\begin{lemm}\label{drw}
Si $W$ est une $B$-paire, et si $\dtilde^r(W) = \{ y \in
\btrig{,r}[1/t] \otimes_{\be} W_e$ tels 
que $\iota_n(y) \in W_{dR}^+$
pour tout $n\geq n(r) \}$, alors :
\begin{enumerate}
\item $\dtilde^r(W)$ est un $\btrig{,r}$-module libre de rang $d$;
\item $\dtilde^r(W) [1/t] = \btrig{,r}[1/t] \otimes_{\be} W_e$;
\item $\dtilde^r(W)$ est stable par $G_K$ et $\phi(\dtilde^r(W)) = 
\dtilde^{pr}(W)$.
\end{enumerate}
\end{lemm}

\begin{proof}
Commençons par remarquer que si $n \geq n(r)$, alors l'image de
$\btrig{,r}[1/t] \otimes_{\be} W_e$ par l'application $\iota_n$
est dense dans $W_{dR}$.
Soit $\dtilde^r_{n(r)}(W) = \{ y \in
\btrig{,r}[1/t] \otimes_{\be} W_e$ tels que $\iota_{n(r)}(y) \in 
W_{dR}^+ \}$; c'est un $\btrig{,r}$-module
libre de rang $d$ (il est engendré par $d$ éléments dont les images forment une base de $W_{dR}^+$). Par ailleurs $\dtilde^r(W)$ est fermé dans 
$\dtilde^r_{n(r)}(W)$ et par la proposition \ref{btpbez},  
$\dtilde^r(W)$ est libre de rang $\leq d$. 

Montrons que si $x \in \btrig{,r}[1/t] \otimes_{\be} W_e$, alors il
existe $m \geq 0$ tel que $t^m x \in \dtilde^r(W)$ ce qui implique que
$\dtilde^r(W)$ est un $\btrig{,r}$-module libre de rang $d$
et que $\dtilde^r(W) [1/t] = \btrig{,r}[1/t] \otimes_{\be} W_e$.
Si $e_1,\dots,e_d$ est une base de $W_e$, alors il existe $m_1$ tel
que pour tout $n \geq n(r)$, l'image par $\iota_n$ des $t^{m_1} e_i$
appartient à $W_{dR}^+$ (si c'est vrai pour un $n$, c'est vrai pour tous 
car $\phi(e_i)=e_i$). Comme $\be = \cup_{j \geq 0} 
t^{-j} (\btp)^{\phi=p^j}$, il existe $m_2$ tel que $t^{m_2} x$ appartient
au $\btrig{,r}$-module engendré par les $e_i$. On peut alors prendre 
$m = m_1 + m_2$.

Ceci démontre les points (1) et (2), et le (3) est une évidence. 
\end{proof}

\begin{lemm}\label{imdr}
Pour tout $n \geq n(r)$, l'image de $\dtilde^r(W)$ par $\iota_n$ contient
une base de $W_{dR}^+$ et si $\dtilde'$ est un
sous-$\btrig{,r}$-module fermé stable par $G_K$ de $\dtilde^r(W)$, dont
l'image par $\iota_n$ contient une base de $W_{dR}^+$ 
pour tout $n \geq n(r)$, alors $\dtilde' = \dtilde^r(W)$.
\end{lemm}

\begin{proof}
Soient $n \geq n(r)$ et $x_1,\dots,x_d$ des éléments de $\btrig{,r}[1/t] \otimes_{\be} W_e$ dont les images par $\iota_n$ forment une base de
$W_{dR}^+$. Si $\ell \geq 0$ est tel que $t^\ell x_i \in \dtilde^r(W)$ pour tout $i$, alors posons $y_i = (t/Q_n)^\ell x_i$. Pour tout $m \geq n(r)$
on a $\iota_m(y_i) \in W_{dR}^+$ et par ailleurs $\iota_n(y_i) = 
\mathrm{inversible} \cdot \iota_n(x_i)$ ce qui fait que l'image de $\dtilde^r(W)$ par $\iota_n$ contient bien une base de $W_{dR}^+$.

Pour montrer l'unicité, on se ramène au cas de rang $1$ en prenant le 
déterminant. Il faut donc montrer que si $x \in \btrig{,r}$ engendre un idéal
stable par $G_K$ et si $\iota_n(x)$ est une unité de $\bdr^+$ pour tout
$n \geq n(r)$, alors $x$ est une unité. Soit $\eta : G_K \to 
(\btrig{,r})^\times = (\btdag{,r})^\times$ 
le cocycle $g \mapsto g(x)/x$. Par la proposition 4.2.1 de
\cite{BC}, l'anneau $\btdag{,r}$ satisfait
les conditions de Tate-Sen et il existe donc une extension finie $L$ de $K$,
un entier $m \geq 0$ et une unité $y \in (\btdag{,r})^\times$ tels que
$h(xy)=xy$ si $h \in H_L$ et $g(xy)/xy \in \phi^{-m}(\bdag{,p^m r}_L)$
si $g \in G_K$. Le lemme 3.2.5 de \cite{BC} montre alors que, quitte
à augmenter $m$, on a $xy \in \phi^{-m}(\bdag{,p^m r}_L)$. Si $L_0 
\subset L$ est une sous-extension non ramifiée de $L$, alors 
$\mathrm{N}_{L/L_0}(\phi^m(xy))$ est un élément de 
$\bdag{,p^m r}_{L_0}$ qui engendre un idéal stable par un
sous-groupe ouvert de $\Gamma_L$. Un raisonnement analogue à celui
du lemme I.3.2 de \cite{Ber3} montre que cet idéal est engendré
par un élément de la forme $\prod_{n \geq n(r)} 
(Q_{n+m}/p)^{\alpha_n}$ et la condition selon laquelle 
$\iota_n(x)$ est une unité de $\bdr^+$ pour tout
$n \geq n(r)$ nous dit que les $\alpha_n$ sont tous nuls, ce qui fait
que $\phi^m(xy)$ est une unité, et donc que $x$ est une unité.
\end{proof}

En particulier, si $s \geq r$, alors l'application naturelle $\btrig{,s} 
\otimes_{\btrig{,r}} \dtilde^r(W) \to \dtilde^s(W)$ est un isomorphisme.

\begin{defi}\label{defdt}
On définit $\dtilde(W) = \btrig{} 
\otimes_{\btrig{,r}} \dtilde^r(W)$ et
si $I$ est un intervalle contenu dans $[r;+\infty[$, 
alors on pose $\dtilde^I(W) =
\bt^I \otimes_{\btrig{,r}} \dtilde^r(W)$. 
\end{defi}

Le lemme \ref{imdr} ci-dessus
montre que cela ne dépend pas du choix de $r \in I$. Remarquons en particulier que si $J \subset I$, alors $\dtilde^J(W) =
\bt^J \otimes_{\bt^I} \dtilde^I(W)$.

\begin{prop}\label{techpbi}
Si $W$ est une $B$-paire, et si $I$ est un intervalle, alors
il existe $j \geq 0$ et une extension finie $L$ de $K$ tels que pour tout $k \geq 0$, il existe un $\bfont^{p^{j+k}I}_L$-module 
$\dfont^{p^{j+k}I}_L$ libre de rang fini vérifiant : 
\begin{enumerate}
\item $\bt^{p^{j+k}I} 
\otimes_{\bfont^{p^{j+k}I}_L} \dfont^{p^{j+k}I}_L =
\dtilde^{p^{j+k}I}(W)$; 
\item $\phi^*(\dfont^{p^{j+k}I}_L) 
= \dfont^{p^{j+k+1}I}_L$; 
\item les images de $\dfont^{p^{j+k}I}_L$
et $\dfont^{p^{j+k+1}I}_L$ dans 
$\dtilde^{p^{j+k}I \cap p^{j+k+1}I}(W)$ engendrent le même
$\bfont^{p^{j+k}I \cap p^{j+k+1}I}_L$-module.
\end{enumerate}
\end{prop}

\begin{proof}
Si l'on se donne une base de $\dtilde^I(W)$, alors l'action de $G_K$ est 
donnée par une application $G_K \to \GL_d(\bt^I)$ et il existe une extension
finie $L$ de $K$ telle que l'image de $G_L$ par 
cette application soit incluse dans
l'ensemble des matrices $M$ vérifiant $V_I(1-M) > c_1 + 2 c_2 + 2 c_3$.
Par la proposition \ref{tsbrs}, $\bt^I$ satisfait les conditions de Tate-Sen
et la proposition 3.2.6 de \cite{BC} nous fournit alors une nouvelle base 
de $\dtilde^I(W)$ et $n \geq 0$ tels 
que l'application $G_K \to \GL_d(\bt^I)$ est triviale
sur $H_L$ et à valeurs dans $\GL_d(\phi^{-n}(\bfont^{p^nI}_L))$.
Si $\dfont^{p^n I}_L$ est le $\bfont^{p^nI}_L$-module engendré par 
$\phi^n$ de cette base, alors $\dtilde^{p^n I}(W) = \bt^{p^n I}
\otimes_{\bfont^{p^nI}_L}  \dfont^{p^n I}_L$. 

Posons $J = p^n I \cap p^{n+1} I$ et $\dfont^{p^{j+k+1} I}_L = \phi^* (\dfont^{p^{j+k} I}_L)$ pour $k \geq 0$. Si $J$ est vide, alors il suffit de prendre $j = n$, les conditions (1) et (2) étant vérifiées et la condition (3) vide. 
Supposons donc que $J$ est non vide; l'unicité dans la méthode de Sen (cf.\ 
la démonstration du (3) de la proposition 3.3.1 de \cite{BC}) ne nous
donne pas la condition (3), mais nous dit qu'il existe $m \geq 0$ tel que :
\[ \phi^{-m}(\bfont^{p^m J}_L) \otimes_{\bfont^{p^nI}_L} 
\dfont^{p^nI}_L = \phi^{-m}(\bfont^{p^m J}_L) 
\otimes_{\bfont^{p^{n+1} I}_L}  \dfont^{p^{n+1}I}_L. \]
En appliquant $\phi^m$ à cette relation, on voit que (3) est satisfaite
en prenant $j=m+n$.
\end{proof}

\begin{prop}\label{pgok}
Si $W=(W_e,W_{dR}^+)$ est une $B$-paire de dimension $d$,
alors il existe un unique $(\phi,\Gamma_K)$-module $\dfont(W)$ 
sur $\brig{}{,K}$ tel que $\btrig{} \otimes_{\brig{}{,K}} 
\dfont(W) = \dtilde(W)$.
\end{prop}

\begin{proof}
Choisissons un intervalle $I$ tel que $I \cap pI$ est non vide. Le (3) de la proposition \ref{techpbi} nous fournit une collection compatible de 
$\bfont^{p^{j+k}I}_L$-modules, et par la définition 2.8.1 et le théorème
2.8.4 de \cite{KK05}, il existe un $\brig{}{,L}$-module $\dfont_L(W)$ 
libre de rang $d$ et tel que $\dfont^{p^{j+k}I}_L  = 
\bfont^{p^{j+k}I}_L \otimes_{\brig{}{,L}} \dfont_L(W)$ pour $k \geq 0$.
La condition (1) implique que $\btrig{} \otimes_{\brig{}{,L}} 
\dfont_L(W) = \dtilde(W)$ et la condition (2) implique que $\dfont_L(W)$
est un $\phi$-module. Par ailleurs $\dfont_L(W)$ est stable sous
l'action de $G_K$ et $H_L$ agit trivialement dessus, puisque c'est vrai pour chaque $\dfont^{p^{j+k}I}_L$. Enfin, le lemme 4.2.5 de \cite{BC}
montre que l'extension $\brig{}{,L} / \brig{}{,K}$ vérifie les conditions 
de la proposition 2.2.1 de descente étale 
de \cite{BC} ce qui fait que l'on a $\dfont_L(W) 
= \brig{}{,L} \otimes_{\brig{}{,K}} \dfont_L(W)^{H_K}$ et on peut donc prendre $\dfont(W)  =  \dfont_L(W)^{H_K}$. Ceci montre l'existence
de $\dfont(W)$.

Passons à l'unicité. Si $D_1$ et $D_2$ satisfont les conclusions de la 
proposition, choisissons des bases de $D_1$ et $D_2$, et
appelons $G_i$ la matrice de $\gamma \in \Gamma_K$
sur $D_i$ et $M$ la matrice de passage d'une base à l'autre, ce qui
fait que $G_1 M=\gamma(M) G_2$. On se donne $r \gg 0$
tel que toutes ces matrices ont leurs coefficients dans $\btrig{,r}$ 
et $I$ un intervalle contenu dans $[r;+\infty[$. La partie {\og unicité \fg}
de la méthode de Sen nous dit qu'il existe $n \geq 0$ tel que 
$M \in \GL_d(\phi^{-n}(\bfont^{p^n I}_K))$. Comme 
$\bfont^{p^n I}_K \cap \btrig{,p^n r} = \brig{,p^n r}{,K}$ par 
le lemme \ref{intri}, on trouve que $M \in \GL_d(\phi^{-n}( 
\brig{,p^n r}{,K}))$. Si $P_1$ et $P_2$ sont les matrices de $\phi$
sur $D_1$ et $D_2$, alors $P_1 M = \phi(M) P_2$ et donc si
$M \in \GL_d(\phi^{-n}(\brig{}{,K}))$, alors 
$M \in \GL_d(\brig{}{,K})$, ce qui fait que $D_1=D_2$.
\end{proof}

Rappelons que  
si $\dfont$ est un $\phi$-module sur $\brig{}{,K}$, alors 
pour $r\gg 0$ on note $\dfont^r$ le $\brig{,r}{,K}$-module fourni 
par le théorème I.3.3 de \cite{Ber5}. En particulier
$\dfont^{pr} = 
\brig{,pr}{,K} \otimes_{\brig{,r}{,K}} \dfont^r$ 
et $\dfont^{pr}$ a une base contenue dans $\phi(\dfont^r)$.
Par exemple, $\dfont^r(W) = \dtilde^r(W) \cap \dfont(W)$
si $r \gg 0$.

\begin{prop}\label{invpg}
Si $\dfont$ est un $(\phi,\Gamma_K)$-module de
rang $d$ sur $\brig{}{,K}$, alors :
\begin{enumerate}
\item $W_e(\dfont) = 
(\btrig{}[1/t] \otimes_{\brig{}{,K}} 
\dfont)^{\phi=1}$ est un $\be$-module 
libre de rang $d$ stable sous l'action de $G_K$; 
\item $W_{dR}^+(\dfont)
=\bdr^+ \otimes^{\iota_n}_{\brig{,r_n}{,K}} 
\dfont^{r_n}$ ne dépend 
pas de $n \gg 0$ et c'est un $\bdr^+$-module 
libre de rang $d$ stable sous l'action de $G_K$;
\item le couple $W(\dfont)=(W_e(\dfont),W_{dR}^+(\dfont))$ est 
une $B$-paire.
\end{enumerate}
\end{prop}

\begin{proof}
Il est clair que $W_e(\dfont)$ est un $\be$-module stable sous l'action de $G_K$. Reste à montrer qu'il est libre de rang $d$. Pour cela, soit
$\oplus M_{a_i,h_i}$ une décomposition de $\btrig{} 
\otimes_{\brig{}{,K}} \dfont$ en $\phi$-modules élémentaires fournie par 
la proposition \ref{dmked}. On a $W_e(\dfont) =  \oplus (M_{a_i,h_i}
[1/t])^{\phi=1}$. On vérifie que l'application $(M_{a,h} [1/t])^{\phi=1} 
\to \btp[1/t]$ qui à $\sum_{i=0}^{h-1} x_i e_i$ associe 
$x_0$ est un isomorphisme 
entre $(M_{a,h} [1/t])^{\phi=1}$ et $(\btp[1/t])^{\phi^h=p^{-a}}$, 
et le (1) suit alors du lemme \ref{bphd}.

Pour montrer le (2), remarquons que par le 
théorème I.3.3 de \cite{Ber5}, on a 
$\phi^*(\dfont^{r_n}) = \dfont^{r_{n+1}}$ pour $n \gg 0$ et
donc le $\bdr^+$-module engendré par $\phi^{-n}(\dfont^{r_n})$ ne
dépend pas de $n \gg 0$. Comme $\dfont^{r_n}$ est libre de rang $d$
et stable par $G_K$, il en est de même pour $W_{dR}^+(\dfont)$.

Montrons maintenant le (3) : $W_{dR}(\dfont) = \bdr \otimes_{\be} 
W_e(\dfont)$ est un $\bdr$-espace 
vectoriel de dimension $d$, réunion croissante des $\bdr$-espaces 
vectoriels $\bdr \otimes_{\be} 
(\btrig{,r_n}[1/t] \otimes_{\brig{,r_n}{,K}} 
\dfont^{r_n})^{\phi=1}$ ce qui fait que si $n \gg 0$, alors :
\[ W_{dR}(\dfont) = \bdr \otimes_{\be} 
(\btrig{,r_n}[1/t] \otimes_{\brig{,r_n}{,K}} 
\dfont^{r_n})^{\phi=1} \subset \bdr 
\otimes^{\iota_n}_{\brig{,r_n}{,K}} 
\dfont^{r_n}
= \bdr \otimes_{\bdr^+} W_{dR}^+(\dfont). \]
En comparant les dimensions, on voit que l'on a en fait égalité et donc
que $W_{dR}^+(\dfont)$ est un réseau de $W_{dR}(\dfont)$.
\end{proof}

\begin{theo}\label{bepgeq}
Les foncteurs $W \mapsto \dfont(W)$ et $\dfont \mapsto W(\dfont)$ 
sont inverses
l'un de l'autre et donnent une équivalence de catégories entre la catégorie
des $B$-paires et la catégorie des 
$(\phi,\Gamma_K)$-modules sur $\brig{}{,K}$.
\end{theo}

\begin{proof}
On vérifie que ces deux foncteurs sont inverses l'un de l'autre en utilisant 
le fait que :
\[ \btrig{}[1/t] \otimes_{\brig{}{,K}} \dfont(W) =
\btrig{}[1/t] \otimes_{\be} W \]
dans la proposition \ref{pgok} et que :
\[ \btrig{}[1/t] \otimes_{\be} W_e(\dfont)
= \btrig{}[1/t] \otimes_{\brig{}{,K}} \dfont \]
dans la proposition \ref{invpg}, puis en identifiant les différents
objets (c'est un exercice instructif que nous laissons au lecteur; il
faut utiliser l'unicité dans la proposition \ref{pgok}).
\end{proof}

\begin{rema}\label{benodrp}
Les modifications de $(\phi,\Gamma_K)$-modules correspondent
à des modifications de $B$-paires. Par exemple : 
\begin{enumerate}
\item Si $W$ et $X$ sont deux $B$-paires, alors $W_e \simeq X_e$ 
si et seulement si $\dfont(W)[1/t] \simeq \dfont(X)[1/t]$;
\item si $W=(W_e,W_{dR}^+)$ est une $B$-paire, 
alors $\dfont((W_e,t^n W_{dR}^+)) = t^n \dfont(W)$.
\end{enumerate}
\end{rema}

\begin{prop}\label{isobp}
Le foncteur $\dfont \mapsto W(\dfont)$ réalise une équivalence de catégories
entre la catégorie des $(\phi,\Gamma_K)$-modules étales et la catégorie
des $B$-paires de la forme $W(V) = (\be \otimes_{\Qp} V, \bdr^+ 
\otimes_{\Qp} V)$ où $V$ est une représentation $p$-adique de $G_K$.
\end{prop}

On retrouve alors (en appliquant le (b) du théorème 6.3.3 de \cite{KK05},
qui nous dit que le foncteur naturel de la catégorie des $(\phi,\Gamma_K)$-modules étales sur $\bdag{}_K$ vers la catégorie des $(\phi,\Gamma_K)$-modules étales sur $\brig{}{,K}$ est une équivalence de catégories) 
le résultat principal de \cite{CC98}, 
c'est-à-dire l'équivalence de catégories entre représentations
$p$-adiques et $(\phi,\Gamma_K)$-modules étales,
avec une démonstration proche de 
celle de \cite{BC}.

Dans le paragraphe \ref{subisoc}, nous reviendrons sur le 
problème de la description {\og explicite \fg} des $B$-paires
dont le $(\phi,\Gamma)$-module associé est isocline.

\subsection{Théorie de Hodge $p$-adique}\label{subthp}
Dans ce paragraphe, nous généralisons  les notions 
habituelles de théorie de Hodge $p$-adique aux $B$-paires.

\begin{defi}\label{adm}
Si $\heartsuit \in \{ \mathrm{cris}, \mathrm{st}, \mathrm{dR} \}$, 
et si $W$ est une $B$-paire,
alors on dit que $W$ est \emph{cristalline}
(ou \emph{semi-stable} ou \emph{de de Rham}) 
si la $\bfont_\heartsuit$-représentation 
$\bfont_\heartsuit \otimes_{\be} W_e$ est triviale. On
pose $\dfont_\heartsuit(W) = 
(\bfont_\heartsuit \otimes_{\be} W_e)^{G_K}$.
\end{defi}

Remarquons que si $V$ est une représentation $p$-adique,
alors bien sûr $V$ est cristalline (ou semi-stable, ou de de Rham, 
ou de Hodge-Tate)
si et seulement si $W(V)$ l'est.

\begin{lemm}\label{dstbe}
Si $D$ est un $(\phi,N)$-module sur $K_0$, 
alors $(\bst \otimes_{K_0} D)^{\phi=1,N=0}$ est
un $\be$-module libre de rang $d = \dim(D)$.
\end{lemm}

\begin{proof}
Si $u \in \bst$ vérifie $N(u)=1$, alors on a $\bst=\bmax[u]$ et
si $x \in (\bst \otimes_{K_0} D)^{N=0}$, alors on peut 
écrire $x=x_0 + u x_1 + \cdots + u^e x_e$ avec $x_i \in
\bmax \otimes_{K_0} D$. Par la proposition 11.7 de \cite{CEV}, 
l'application $x \mapsto x_0$ de $(\bst \otimes_{K_0} D)^{N=0}$ 
dans $\bmax \otimes_{K_0} D$ est un isomorphisme, et on se ramène
donc à montrer que si $D$ est un $\phi$-module sur $K_0$, 
alors $(\bmax \otimes_{K_0} D)^{\phi=1}$ est
un $\be$-module libre de rang $d$. Par le théorème
de Dieudonné-Manin, $\Qpnrhat \otimes_{K_0} D$ se décompose
en somme directe de $\phi$-modules élémentaire $M_{a_i,h_i}$.
L'application $x=\sum_{j=0}^{h-1} x_j e_j \mapsto x_0$ 
de $(\bmax \otimes_{\Qpnrhat} M_{a,h})^{\phi=1}$ dans
$\bmax^{\phi^h=p^{-a}}$ est un isomorphisme, et le lemme
résulte alors du lemme \ref{bphd}.
\end{proof}

\begin{prop}\label{usuthp}
Si $W$ est une $B$-paire, alors :
\begin{enumerate}
\item $\dst(W)$ est un $(\phi,N)$-module sur $K_0$ et 
$\dcris(W)=\dst(W)^{N=0}$;
\item $\ddr(W)$ est un $K$-espace vectoriel filtré 
avec $\fil^i(\ddr(W)) = \ddr(W) \cap t^i W_{dR}^+$
et l'application
naturelle $K \otimes_{K_0} \dst(W) \to \ddr(W)$ est injective;
\item si $W$ est semi-stable, alors $W_e = (\bst \otimes_{K_0} 
\dst(W))^{\phi=1,N=0}$ et $W_{dR}^+ = \fil^0(\bdr \otimes_K
\ddr(W))$;
\item si $D$ est un $(\phi,N)$-modules filtré, et si
$W_e(D) = (\bst \otimes_{K_0} D)^{\phi=1,N=0}$ 
et $W_{dR}^+(D) = \fil^0(\bdr \otimes_K D)$, alors
$W(D)=(W_e(D),W_{dR}^+(D))$ est une $B$-paire semi-stable.
\end{enumerate}
\end{prop}

\begin{proof}
Exercice (ce n'est pas différent du cas des représentations 
$p$-adiques semi-stables).
\end{proof}

A la lumière du théorème \ref{bepgeq}, l'étude des $B$-paires de de Rham 
revient à l'étude des $(\phi,\Gamma_K)$-modules
{\og localement triviaux \fg} au sens de \cite{Ber5}, ce qui est fait
en détail dans \cite{Ber5}.

Si $D$ est un $(\phi,N)$-module filtré, on note $\calM(D)$
le $(\phi,\Gamma_K)$-module construit dans \cite[\S II.2]{Ber5}.
Rappelons que :
\[ \calM(D) = \{ y \in (\brig{}{,K}[\log(X),1/t] \otimes_{K_0} D)^{N=0}
\ \mid\  \iota_n(y) \in \fil^0(K_n(\!(t)\!) \otimes_K D_K) \ \forall n \gg 0\}. \]

\begin{prop}\label{eqsst}
Les foncteurs $W \mapsto \dst(W)$ 
et $D \mapsto W(D)$ sont inverses l'un de l'autre et donnent 
une équivalence de catégories entre la catégorie des $B$-paires
semi-stables et la catégorie des $(\phi,N)$-modules filtrés.

Si $W$ est une $B$-paire, alors $\dfont(W) = \calM(\dst(W))$ et donc
si $D$ est un $(\phi,N)$-module filtré, alors $W(D) = W(\calM(D))$.
\end{prop}

\begin{proof}
Ces affirmations ne présentent aucune difficulté.
\end{proof}

\begin{theo}\label{monetadm}
Le théorème de monodromie $p$-adique et le théorème {\og faiblement
admissible implique admissible \fg} sont vrais. De fait,
\begin{enumerate}
\item toute $B$-paire de de Rham est potentiellement semi-stable; 
\item le foncteur $W \mapsto \dst(W)$ 
réalise une équivalence de catégories
entre la catégorie des objets de la 
forme $W(V)$ où $V$ est une représentation
$p$-adique et la catégorie des 
$(\phi,N)$-modules filtrés admissibles.
\end{enumerate}
\end{theo}

\begin{proof}
Montrons tout d'abord le (1).
Soit $W$ une $B$-paire et $\dfont(W)$ le $(\phi,\Gamma_K)$-module
associé. Si $W$ est de de Rham, alors pour $n \gg 0$, 
la $\bdr$-représentation $\bdr \otimes^{\iota_n}_{\brig{,r_n}{,K}}
\dfont^{r_n}(W)$ est égale à $\bdr \otimes_{\be} W_e$ et donc 
triviale, ce qui fait par le théorème 3.9 de \cite{F04} que
le $K_\infty(\!(t)\!)$-module à connexion $K_\infty(\!(t)\!) 
\otimes^{\iota_n}_{\brig{,r_n}{,K}}
\dfont^{r_n}(W)$ est trivial. Etant donnée
la définition III.1.2 de \cite{Ber5}, on est en mesure d'en appliquer 
le théorème A de \cite{Ber5}, qui nous dit qu'il existe une
extension finie $L$ de $K$ et un $(\phi,N)$-module filtré $D$ sur $L$
tels que $\dfont(W_{\mid L}) = \calM(D)$ et donc que $W_{\mid L}$
est semi-stable. Ceci montre le (1). Le (2) suit du théorème
A de \cite{CF} ou bien (si l'on préfère 
passer par les $(\phi,\Gamma_K)$-modules)
du théorème B de \cite{Ber5}.
\end{proof}

Rappelons que si $U$ est une $\Cp$-représentation de $G_K$, alors
la réunion $U^{H_K}_{\mathrm{fini}}$ 
des sous-$K_\infty$-espaces vectoriels de dimension finie stables
par $\Gamma_K$ de $U^{H_K}$ a la propriété que
l'application  $\Cp \otimes_{K_\infty}
U^{H_K}_{\mathrm{fini}} \to U$ est un isomorphisme (cf.\ \cite{Sn80}).
L'espace $U^{H_K}_{\mathrm{fini}}$ est muni de l'application 
$K_\infty$-linéaire $\nabla_U = \log(\gamma) / \log_p(\chi(\gamma))$
avec $\gamma \in \Gamma_K \setminus \{1\}$ suffisamment proche de $1$.

\begin{defi}\label{senht}
Si $W$ est une $B$-paire, alors $W_{dR}^+ / t W_{dR}^+$ est une
$\Cp$-représentation de $G_K$ et on 
pose $\dsen(W) = (W_{dR}^+ / t W_{dR}^+)^{H_K}_{\mathrm{fini}}$.
On pose $\Theta_{\mathrm{Sen}} = \nabla_W$ et on 
dit que $W$ est \emph{de Hodge-Tate} si 
$\Theta_{\mathrm{Sen}}$ est diagonalisable
à valeurs propres appartenant à $\ZZ$. Ces entiers sont 
les \emph{poids de Hodge-Tate} de $W$.
\end{defi}

\section{Les $(\phi,\Gamma)$-modules}\label{appber}
Etant donné le théorème \ref{bepgeq}, l'étude des $B$-paires revient à
l'étude des $(\phi,\Gamma)$-modules. Dans ce chapitre, nous montrons
plusieurs résultats sur les $(\phi,\Gamma)$-modules : modification,
classification des objets isoclines, classification des objets de hauteur finie.

\subsection{Modification de $(\phi,\Gamma)$-modules}\label{submodpg}
Si $\dfont$ est un $\phi$-module 
sur $\brig{}{,K}$ et si $n \geq n(r)$ avec 
$r \gg 0$, rappelons que $\dfont^r / Q_n$ est un $K_n$-espace
vectoriel de dimension $\rg(D)$. Soit
$M = \{M_n\}_{n \geq n(r)}$ une famille $\phi$-compatible de 
sous-$K_n$-espaces vectoriels de $\dfont^r / Q_n$, c'est-à-dire que pour
tout $n \geq n(r)$, $M_{n+1}$ est engendré par les $\phi(y)$ où $y 
\in \dfont^r$ est tel que son image dans $\dfont^r / Q_n$ appartient à $M_n$.
En d'autres termes, on un isomorphisme 
$K_{n+1} \otimes_{K_n} (\dfont^r / Q_n) \to \dfont^{pr} / Q_{n+1}$
obtenu en quotientant l'isomorphisme 
$\phi^*(\dfont^r) \simeq
\dfont^{pr}$ par $\phi(Q_n) = Q_{n+1}$, et qui est donc 
donné par $\alpha \otimes \overline{m} \mapsto \alpha \otimes
\overline{\phi(m)}$ et on demande que $M_{n+1}$ soit l'image
de $K_{n+1} \otimes_{K_n} M_n$ par cet isomorphisme (voir 
\cite[\S II.1]{Ber5} pour une condition analogue). 

\begin{defi}\label{donmod}
Une telle famille $M = \{M_n\}_{n \geq n(r)}$
de sous-espaces vectoriels de $\dfont^r / Q_n$ est appelée une
\emph{donnée de modification} de $\dfont$. On définit alors
$\dfont[M] = \{ y \in \dfont$ dont l'image dans $\dfont^r/Q_n$ appartient à 
$M_n \}$.
\end{defi} 

\begin{prop}\label{hammer}
Si $\dfont$ est un $\phi$-module de rang $d$ 
sur $\brig{}{,K}$ et si 
$M$ est une donnée de modification, 
alors $\dfont[M]$ est un $\phi$-module et de plus :
\begin{enumerate}
\item si $M \subset N$, alors $\dfont[M] \subset \dfont[N]$ et si 
$\dfont[M]=\dfont[N]$, alors $M=N$;
\item $\dfont[M]=\dfont$ si $M_n = \dfont^r/Q_n$ pour tout $n$, 
et $\dfont[0]=t \cdot \dfont$;
\item $\deg(\dfont[M]) = \deg(\dfont)+d-\dim(M)$;
\item si $\dfont$ est un $(\phi,\Gamma_K)$-module et si
$M$ est stable par $\Gamma_K$, alors $\dfont[M]$ 
est un $(\phi,\Gamma_K)$-module.
\end{enumerate}
\end{prop}

\begin{proof}
Rappelons que $t=\log(1+X)=X \cdot \prod_{n \geq 0} Q_n/p$ 
ce qui fait que $t\cdot \dfont \subset \dfont[M]$ quel que soit $M$.
La définition de $\dfont[M]$ implique que c'est un sous-module 
fermé de $\dfont$, et comme il contient $t\cdot \dfont$, il est libre de
rang $d$. Le fait que c'est un $\phi$-module résulte du fait
que la famille des $M_n$ est $\phi$-compatible. 

Le fait que si $M \subset N$, alors $\dfont[M] \subset \dfont[N]$ est
évident. Remarquons que si l'image de $x \in \dfont$ appartient à $M_n$,
alors $x \cdot t/Q_n  \in \dfont[M]$ et son image dans $\dfont^r/Q_n$ est 
(un multiple de) $x$  ce qui fait que l'image de $\dfont[M]$ dans $\dfont^r/Q_n$
est $M_n$ et en particulier que si $\dfont[M]=\dfont[N]$, alors $M=N$. 
Ceci montre le (1). Le (2) est une évidence.

Pour montrer le (3), prenons un drapeau $\{0\} 
= M^{(0)} \subset M^{(1)} \subset \cdots \subset M^{(d)}$
tel que $M^{(\dim(M))} = M$. On en déduit une suite 
de $\phi$-modules $t\cdot \dfont = \dfont^{(0)} \subset \cdots \subset \dfont^{(d)} 
=\dfont$ et donc une suite de $\phi$-modules de rang $1$ : 
\[ \det(t\cdot \dfont) = \det(\dfont^{(0)}) \subset 
\cdots \subset \det(\dfont^{(d)}) 
= \det(\dfont). \]
Leurs pentes forment une suite strictement décroissante 
(si on a égalité des pentes, on a égalité des modules par le lemme
\ref{above}) de $d+1$ nombres entiers (ils sont entiers 
par le lemme \ref{unis})
dont le premier est $d + \deg(\dfont)$ et le dernier 
est $\deg(\dfont)$, ce qui fait que $\deg(\dfont^{(m)})$ est
forcément $\deg(\dfont)+d-\dim(M)$. Enfin, le (4) est une évidence.
\end{proof}

\begin{prop}\label{modext}
Si $M$ est une donnée de modification de $\dfont$ et si l'on a une suite
exacte $0 \to \dfont_1 \to \dfont \to \dfont_2 \to 0$, alors :
\begin{enumerate}
\item $\{ M_n \cap \dfont_1^r/Q_n \}$ est une donnée 
de modification pour $\dfont_1$ et $\{ \im(M_n \to
\dfont_2^r/Q_n) \}$ est une donnée de modification pour $\dfont_2$;
\item on a une suite exacte $0 \to \dfont_1[M] \to \dfont[M] \to \dfont_2[M] \to 0$.
\end{enumerate}
\end{prop}

\begin{proof}
Le (1) suit du fait que les $M_n$ sont $\phi$-compatibles et
que les $\dfont_i^r / Q_n$ le sont aussi. Passons au (2); on a clairement
une suite exacte $0 \to \dfont_1[M] \to \dfont[M] \to \dfont_2[M]$.
Par le (3) de la proposition \ref{hammer}, on a :
\begin{enumerate}
\item 
$\deg(\dfont_1[M]) = \deg(\dfont_1) + \dim(\dfont_1) - \dim(M_n \cap \dfont_1^r/Q_n)$; 
\item  $\deg(\dfont[M]) = \deg(\dfont) + \dim(\dfont) - \dim(M_n)$;
\item  $\deg(\dfont_2[M]) = \deg(\dfont_2) + \dim(\dfont_2) - 
\dim(\im(M_n \to \dfont_2^r/Q_n))$.
\end{enumerate}
On en déduit que $\deg(\dfont[M]) = \deg(\dfont_1[M]) + \deg(\dfont_2[M])$, et le corollaire
\ref{ptext} implique alors que la suite 
$0 \to \dfont_1[M] \to \dfont[M] \to \dfont_2[M] \to 0$ est exacte.
\end{proof}

\begin{rema}\label{compos}
Si $M$ et $N$ sont deux données de modification telles que
$M \cap N = 0$, alors on un isomorphisme $\dfont[M][N] = \dfont[M \oplus N]$ 
(il y a une inclusion évidente; comparer les degrés).
\end{rema}

\begin{theo}\label{subet}
Si $\dfont$ est un $\phi$-module sur $\brig{}{,K}$, alors il existe 
un $\phi$-module étale $\dfont' \subset \dfont[1/t]$ tel que $\dfont'[1/t]  = \dfont[1/t]$.
\end{theo}

\begin{proof}
Il suffit de montrer qu'il existe un nombre fini de modifications 
successives (au sens de la proposition \ref{hammer})
de $\dfont$ dont le résultat est isocline de pente entière
(si $\dfont$ est isocline de pente $s$, alors $t \cdot \dfont$
est isocline de pente $s+1$). 
Remarquons
que si on a une suite exacte $0 \to \dfont_1 \to \dfont \to \dfont_2 \to 0$, alors
il est toujours possible de modifier $\dfont_1$ ou $\dfont_2$ sans toucher à l'autre,
en choisissant $M$ convenablement. 

Etape 1 : si $\dfont$ est isocline, alors $\dfont$ se modifie en une extension
successive de $\phi$-modules de pentes entières. En effet, si $\dfont$ est 
isocline, modifions le par un $M$ de codimension $1$ ce qui augmente 
le degré de $\dfont$ de $1$. Si $\dfont[M]$ est isocline de pente non entière, alors
on répète cette opération. Le résultat est donc, après un nombre fini 
de modifications, que $\dfont$ devient isocline de pente entière (et on a terminé
l'étape 1) ou bien que $\dfont$ se casse en deux morceaux. Dans ce dernier cas, 
on a terminé par récurrence sur la dimension de $\dfont$.

Etape 2 : si $\dfont$ est extension
successive de $\phi$-modules de pentes entières, alors il se modifie en
un $\phi$-module isocline de pente entière. On se ramène au cas d'une
extension de deux $\phi$-modules de pentes entières  
$0 \to \dfont_1 \to \dfont \to \dfont_2 \to 0$ avec $\dfont_i$ isocline de pente $s_i$. 
Dans ce cas :
\begin{enumerate}
\item si $s_2 = s_1$, alors on a terminé par le lemme \ref{extetal};
\item si $s_2 < s_1$, alors la modification est facile à écrire : on remplace
$\dfont$ par $\dfont'$, l'ensemble des $y \in \dfont$ dont l'image dans $\dfont_2$
appartient à $t^{s_1-s_2} \dfont_2$, 
ce qui fait que l'on a une suite exacte
$0 \to \dfont_1 \to \dfont' \to t^{s_1-s_2} 
\dfont_2 \to 0$ et on a terminé par le lemme \ref{extetal};
\item si $s_1 < s_2$, alors on modifie $\dfont$ par une donnée nulle sur
$\dfont_1$ et surjective sur $\dfont_2$ ce qui augmente $s_1$ de $1$ et on
itère cette opération $s_2-s_1$ fois.
\end{enumerate}
\end{proof}

\begin{coro}\label{fontconj}
Si $W_e$ est une $\be$-représentation de $G_K$, alors il existe une
représentation $p$-adique $V$ de $H_K$
telle que la restriction de $W_e$ à $H_K$ est isomorphe à 
$\be \otimes_{\Qp} V$. 
\end{coro}

\begin{proof}
Soit $W_{dR}^+$ un $\bdr^+$-réseau stable par $G_K$ 
de $\bdr \otimes_{\be} W_e$ et $\dfont = 
\dfont(W)$ le $(\phi,\Gamma_K)$-module
associé à la $B$-paire $W=(W_e,W_{dR}^+)$ par le théorème \ref{bepgeq}.
Par le théorème \ref{subet}, il existe un  
$\phi$-module étale
$\dfont' \subset \dfont[1/t]$ 
tel que $\dfont'[1/t]  =  \dfont[1/t]$.  Comme $\dfont'$ est étale,
le $\Qp$-espace vectoriel $V=(\btrig{} \otimes_{\brig{}{,K}} \dfont')^{\phi=1}$
est une représentation de $H_K$ dont la dimension est le rang de $\dfont'$.
On a alors :
\begin{align*} 
W_e & = (\btrig{}[1/t] \otimes_{\brig{}{,K}} \dfont)^{\phi=1} \\
& = (\btrig{}[1/t] \otimes_{\brig{}{,K}} \dfont')^{\phi=1} \\
& = \be \otimes_{\Qp} V, 
\end{align*}
et le corollaire est démontré.
\end{proof}

\begin{rema}\label{bvunik}
\begin{enumerate}
\item Dans le corollaire \ref{fontconj} ci-dessus, on est loin d'avoir unicité. Par
exemple, si $D$ est un $(\phi,N)$-module sur $K_0$ et si $W_e = 
(\bst \otimes_{K_0} D)^{\phi=1,N=0}$, alors $W_e = \be 
\otimes_{\Qp} V$ pour toute représentation $V$ qui s'obtient comme
$V_{\mathrm{st}}(D)$ à partir d'une 
filtration admissible sur $K \otimes_{K_0} D$;
\item Dans la première version de cet article, j'affirmais que si $\dfont$ est un 
$(\phi,\Gamma_K)$-module sur $\brig{}{,K}$, alors il existe une extension 
finie $L$ de $K$ et un $(\phi,\Gamma_L)$-module étale $\dfont'$ sur 
$\brig{}{,L}$ tel que $\dfont' \subset \brig{}{,L}[1/t] \otimes_{\brig{}{,K}} 
\dfont$ et $\brig{}{,L}[1/t] \otimes_{\brig{}{,K}}  
\dfont' = \brig{}{,L}[1/t] \otimes_{\brig{}{,K}} \dfont$. 
Ceci est incorrect, et Kiran Kedlaya et Ruochuan Liu ont construit le 
contre-exemple suivant : soient $p \neq 2$, 
$K=\Qp$ et $\dfont$ une extension non-triviale de 
$\brig{}{,K}$ par $\brig{}{,K}(p^{-2})$ (ce qui veut dire que l'action de $\Gamma_K$
est inchangée et que $\phi$ est multiplié par $p^{-2}$). En utilisant le théorème de
filtration par les pentes de Kedlaya et les calculs de cohomologie de 
$(\phi,\Gamma)$-modules de \cite{PCtr} et de \cite{L07}, 
on peut montrer qu'une telle extension non-triviale existe et que quelque soit 
l'extension finie $L$ de $K$, $\brig{}{,L}  \otimes_{\brig{}{,K}}\dfont$ ne contient
pas de sous-$(\phi,\Gamma_L)$-module isocline de rang $2$ et donc qu'un 
$\dfont'$ comme ci-dessus n'existe pas.
\end{enumerate}
\end{rema}

\subsection{Classification des objets isoclines}\label{subisoc}
Le théorème \ref{bepgeq} nous donne une filtration 
sur les $B$-paires, et il est intéressant de décrire les 
$B$-paires correspondant
aux $(\phi,\Gamma_K)$-modules isoclines. Le cas étale 
relève de la proposition \ref{isobp}.

\begin{defi}\label{rephd}
Si $h \geq 1$ et $a \in \ZZ$ sont premiers entre eux, 
soit $\rep(a,h)$ la catégorie dont 
les objets sont les $\Qph$-espaces vectoriels $V_{a,h}$ 
de dimension finie,
munis d'une action 
semi-linéaire de $G_K$ et d'un Frobenius 
lui aussi semi-linéaire $\phi : V_{a,h} \to V_{a,h}$ qui
commute à $G_K$ et qui vérifie $\phi^h=p^a$. Les morphismes sont
ceux que l'on imagine. 
\end{defi}

\begin{rema}\label{remhd}
\begin{enumerate}
\item Si $V_{a,h} \in \rep(a,h)$, alors $\dim_{\Qph}(V_{a,h})$
est divisible par $h$ (si $e$ est cette dimension, alors
$\phi^h = p^{ae}$ sur $\det(V_{a,h})$);
\item si $h=1$ et $a=0$, alors on retrouve simplement la catégorie
des représentations $p$-adiques de $G_K$; 
\item si $e \in \ZZ$, alors $\rep(a,h)$ et $\rep(a+eh,h)$ sont équivalentes
de manière évidente;
\item si $D$ est un $\phi$-module
isocline de pente $a/h$ sur $K_0$, alors le théorème de Dieudonné-Manin
implique que $V_{a,h} 
= (\Qpnrhat \otimes_{K_0} D)^{\phi^h=p^a}$ est un
objet de $\rep(a,h)$ dont la dimension en tant 
que $\Qph$-espace vectoriel est $\dim_{K_0}(D)$. Dans ce cas, 
l'action de $I_K \subset G_K$ sur $V_{a,h}$ est d'ailleurs triviale.
\end{enumerate}
\end{rema}

Si $V_{a,h} \in \rep(a,h)$, alors on pose $W_e(V_{a,h}) 
= (\bmax \otimes_{\Qph} V_{a,h})^{\phi=1}$ et 
$W_{dR}^+(V_{a,h}) = \bdr^+ \otimes_{\Qph} V_{a,h}$.

\begin{theo}\label{objiso}
Si $V_{a,h} \in \rep(a,h)$, alors $W(V_{a,h}) =(W_e(V_{a,h}),
W_{dR}^+(V_{a,h}))$ est une $B$-paire et le foncteur $V_{a,h}
\mapsto W(V_{a,h})$ définit une équivalence de catégories entre 
$\rep(a,h)$ et la catégorie des $B$-paires $W$ telles que $\dfont(W)$
est isocline de pente $a/h$.
\end{theo}

\begin{proof}
Il est clair que $W = W(V_{a,h})$ est une $B$-paire. Par ailleurs, 
la construction de $\dfont(W)$
fournie par le lemme \ref{drw} et 
la proposition \ref{pgok} montre que l'on a $\btrig{}
\otimes_{\brig{}{,K}} \dfont(W) = \btrig{} 
\otimes_{\Qph} V_{a,h}$ ce qui fait que 
$\dfont(W)$ est isocline de pente $a/h$. 

Réciproquement, si $\dfont$ est un $(\phi,\Gamma_K)$-module
isocline de pente $a/h$, alors par la proposition \ref{dmked}, on
a une décomposition $\btrig{} \otimes_{\brig{}{,K}} \dfont = 
\oplus_{i=1}^k M_{a,h} = \oplus_{i=1}^k \oplus_{j=0}^{h-1}
\btrig{} \phi^j(e_i)$ où $e_i,\phi(e_i),\dots,\phi^{h-1}(e_i)$ est une
base de la $i$-ième copie de $M_{a,h}$. On voit alors que $\sum_{i=1}^k 
\sum_{j=0}^{h-1} \lambda_{ij} \phi^j(e_i) \in 
(\btrig{} \otimes_{\brig{}{,K}} \dfont)^{\phi^h = p^a}$ si et seulement
si on a  $\phi^h(\lambda_{ij}) = \lambda_{ij}$ pour tous $i,j$ et comme
$(\btrig{})^{\phi^h=1} = \Qph$, on trouve que
$V_{a,h} = (\btrig{} \otimes_{\brig{}{,K}} \dfont)^{\phi^h=p^a}$
est un $\Qph$-espace vectoriel de dimension $\rg(\dfont)$
qui hérite d'une action de $G_K$ et d'un Frobenius 
tel que $\phi^h = p^a$. Si $W$
est une $B$-paire telle que $\dfont(W)$
est isocline de pente $a/h$, alors on lui associe l'espace $V_{a,h}$ 
construit à partir de $\dfont(W)$. Le lecteur vérifiera qu'on a 
ainsi défini
un foncteur inverse de $V_{a,h} \mapsto W(V_{a,h})$.
\end{proof}

Il sera intéressant de calculer les extensions d'objets isoclines; les
extensions de $(\phi,\Gamma)$-modules sont étudiées dans \cite{L07}.

\subsection{Les $(\phi,\Gamma)$-modules de hauteur finie}\label{subhf}

Comme l'anneau $\brigplus{,K} = \brig{}{,K} \cap \btp$ 
n'a de bonnes propriétés que si
$K \subset K_0(\mu_{p^\infty})$, on suppose que cette condition
est vérifiée dans tout ce paragraphe. Dans ce cas, $\brigplus{,K}$
s'identifie à l'ensemble des séries $f(X)=\sum_{k \geq 0} f_k X^k$
qui convergent sur le disque unité ouvert, 
ce qui en fait un anneau de Bézout,
et si on pose $\bplus_K = 
\brigplus{,K} \cap \bdag{}_K$, alors $\bplus_K = 
K_0 \otimes_{\OO_{K_0}} 
\OO_{K_0}[\![X]\!]$.

\begin{defi}\label{defhf}
On dit qu'un $(\phi,\Gamma)$-module $\dfont$ sur 
$\brig{}{,K}$ est \emph{de hauteur finie}
s'il existe un $(\phi,\Gamma)$-module $\dfont^+$ sur $\brigplus{,K}$ tel
que $\dfont  = \brig{}{,K} \otimes_{\brigplus{,K}} \dfont^+$ et
on dit qu'une $B$-paire $W$ est de hauteur finie si $\dfont(W)$ l'est.
\end{defi}

Remarquons qu'un $(\phi,\Gamma)$-module sur $\brigplus{,K}$
est un $\brigplus{,K}$-module $\dfont^+$ stable par $\phi$ et 
$\Gamma$ tel que $\det \phi$ est inversible dans $\brig{}{,K}$
(et non dans $\brigplus{,K}$ ce qui serait trop restrictif).

\begin{lemm}\label{chkhf}
Si $K \subset K_0(\mu_{p^\infty})$, alors la définition ci-dessus est compatible avec la définition habituelle quand $W=W(V)$.
\end{lemm}

\begin{proof}
Si $W=W(V)$ avec $V$ de hauteur finie au sens habituel (cf.\ \cite{CH}), 
alors il est évident
que $W$ est de hauteur finie. Montrons donc la réciproque. Par hypothèse,
il existe une base de $\dfont(V)$ dans laquelle $\mat(\phi) = P^+ \in 
\mathrm{M}_d(\brigplus{,K})$ et $\mat(\gamma) = G^+ \in 
\mathrm{GL}_d(\brigplus{,K})$ pour $\gamma \in \Gamma_K$ 
(puisque $V$ est de hauteur finie) ainsi 
qu'une base de $\ddag{}(V)$ 
dans laquelle $\mat(\phi) = P^\dagger \in \GL_d(\bdag{}_K)$ et
$\mat(\gamma) = G^\dagger \in \GL_d(\bdag{}_K)$ pour $\gamma \in \Gamma_K$.
 
Soit $M$ la matrice de passage d'une base à l'autre.
La proposition 6.5 de \cite{KK04} montre que l'on peut 
écrire $M=M^+ \cdot M^\dagger$ avec $M^+ \in 
\GL_d(\brigplus{,K})$ et $M^\dagger \in \GL_d(\bdag{}_K)$.
Dans la base de $\ddag{}(V)$ obtenue en appliquant $M^\dagger$ à celle
de $\ddag{}(V)$, on a :
\begin{align*}
\mat(\phi)&=\phi(M^+)P^+(M^+)^{-1}=
\phi(M^\dagger)^{-1}P^\dagger(M^\dagger), \\
\mat(\gamma)&=\gamma(M^+)G^+(M^+)^{-1}=
\gamma(M^\dagger)^{-1}G^\dagger(M^\dagger), 
\end{align*}
ce qui fait que ces matrices sont 
à coefficients dans $\bdag{}_K \cap \brigplus{,K} = \bplus_K$. 
\end{proof}

\begin{lemm}\label{detpgh}
Si $\dfont^+$ est un $(\phi,\Gamma)$-module sur $\brigplus{,K}$,
alors il existe des entiers $\alpha_0,\dots,\alpha_m$ tels que
l'idéal engendré par $\det(\phi)$ est engendré par
$X^{\alpha_0} Q_1^{\alpha_1} \cdots 
Q_m^{\alpha_m}$ et si $K$ est une extension finie de $\Qp$, alors
$\alpha_0 = 0$.
\end{lemm}

\begin{proof}
Ce lemme se trouve dans la partie B.1.6 de \cite{W96} mais
nous en donnons une démonstration pour la commodité du lecteur.
Comme l'action de $\phi$ commute à $\Gamma_K$, 
l'idéal engendré par $\det(\phi)$ est stable par $\Gamma_K$ 
et la première partie du
lemme suit alors du lemme I.3.2 de \cite{Ber3} et du fait 
que si cet idéal est inversible dans $\brig{}{,K}$, alors
les $\alpha_j$ sont presque tous nuls. 

Si $\delta = \det(\phi)$ et $g = \det(\gamma)$ pour $\gamma
\in \Gamma_K$, alors $\gamma(\delta)/\delta = \phi(g)/g$ et
en réduisant modulo $X$, on trouve que 
$\chi(\gamma)^{\alpha_0}$ s'écrit $\phi(g_0)/g_0$
avec $g_0 \in K_0$. Si
$K$ est une extension finie de $\Qp$, alors $\phi(g_0)/g_0$
est de norme $1$ et donc
$\chi(\gamma)^{\alpha_0 [K_0 : \Qp]} = 1$ 
ce qui fait que $\alpha_0 = 0$.
\end{proof}

Le lecteur vicieux montrera que si le corps résiduel de $K$ est
algébriquement clos, alors on peut effectivement avoir $\alpha_0 \neq 0$.
 
\begin{defi}\label{pgfil}
Un \emph{$\phi$-module filtré sur $K_0$ avec 
action de $\Gamma_K$} est un 
$\phi$-module $D$ sur $K_0$ muni d'une action de $\Gamma_K$
commutant à $\phi$ et d'une filtration 
(décroissante, exhaustive et séparée)
stable par $\Gamma_K$ sur
$D_\infty = K_\infty \otimes_{K_0} D$.
\end{defi}

Il existe alors un entier $n \geq 0$ tel que la filtration de $D_\infty$
est définie sur $K_n$, c'est-à-dire que si l'on pose $D_n =
K_n \otimes_{K_0} D$, alors
$\fil^i D_\infty = K_\infty
\otimes_{K_n} \fil^i D_n$ 
pour tout $i \in \ZZ$
et on appelle $n(D)$
le plus petit entier ayant cette propriété.

\begin{prop}\label{pgfilhf}
Si $D$ est un $\phi$-module filtré sur $K_0$ avec 
action de $\Gamma_K$, alors la
$B$-paire $W=((\bmax \otimes_{K_0} D)^{\phi=1}, 
\bdr^+ \otimes_{K_\infty} D_\infty)$ est de hauteur finie.
\end{prop}

\begin{proof}
Soit $h \geq 0$ tel que $\fil^h D_\infty = 0$ 
et $\fil^{-h} D_\infty = D_\infty$. 
Si $n \geq 1$, on écrit $\phi^{-n}$ pour $\iota_n$. 
Posons alors comme dans \cite[\S 3.1]{PC03} et \cite{Ki06} :
\[ \calM^+(D) = \{ y \in t^{-h} \brigplus{,K} \otimes_{K_0} D
\mid \phi^n(y) \in \fil^0(K_\infty(\!(t)\!) \otimes_{K_\infty} D_\infty) 
\ \forall n \in \ZZ \}. \]
C'est un sous-$\brigplus{,K}$-module fermé de 
$t^{-h} \brigplus{,K} \otimes_{K_0} D$ qui contient
$t^h \brigplus{,K} \otimes_{K_0} D$ et qui est donc libre de
rang $d=\dim(D)$. Il est de plus manifestement stable par $\phi$ 
et $\Gamma_K$.

Si $n \geq n(D)$, alors :
\[ K_n[\![t]\!] \otimes^{\iota_n}_{\brigplus{,K}} 
t^{-h} \brigplus{,K} \otimes_{K_0} D = t^{-h} 
K_n[\![t]\!] \otimes_{K_0} D \supset 
\fil^0(K_n(\!(t)\!) \otimes_{K_n} D_n). \]
Si $y_n \in \fil^0(K_n(\!(t)\!) \otimes_{K_n} D_n)$
et $w \geq 0$, il existe donc 
$y \in t^{-h} \brigplus{,K} \otimes_{K_0} D$
donc l'image par $\iota_n$ vérifie $\iota_n(y)-y_n \in
\fil^w(K_n(\!(t)\!) \otimes_{K_n} D_n)$. 
L'élément $z = y \cdot (t/Q_n)^{2h}$ a alors la propriété
que $\iota_n(z)$ est un multiple de $\iota_n(y)$ et que
pour tout $m \in \ZZ$ différent de $-n$, on a $\phi^m(z)
\in \fil^0(K_\infty(\!(t)\!) \otimes_{K_\infty} D_\infty)$ ce qui fait que
$z \in \calM^+(D)$. On en déduit que si $n \geq n(D)$, alors
l'application : 
\[ K_n[\![t]\!] \otimes^{\iota_n}_{\brigplus{,K}} 
\calM^+(D) \to  \fil^0(K_n(\!(t)\!) \otimes_{K_n} D_n) \]
est un isomorphisme. En particulier, l'application naturelle
$\btrig{} \otimes_{\brigplus{,K}} \calM^+(D)
\to \dtilde(W)$ est un isomorphisme et donc, 
en utilisant l'unicité dans la proposition \ref{pgok}, 
on trouve que $\dfont(W) =
\brig{}{,K} \otimes_{\brigplus{,K}} \calM^+(D)$ ce qui fait
que $W$ est de hauteur finie. Remarquons que dans les notations
de \cite{Ber5}, on a $\calM(D) = \brig{}{,K} 
 \otimes_{\brigplus{,K}} \calM^+(D)$.
\end{proof}

\begin{rema}\label{wthf}
Le fait de ne considérer que les 
$y \in t^{-h} \brigplus{,K} \otimes_{K_0} D$ tels que l'on a 
$\phi^n(y) \in \fil^0(K_\infty(\!(t)\!) \otimes_{K_\infty} D_\infty)$
pour tout $n \in \ZZ$ n'est pas restrictif; le lecteur pourra montrer que :
\[ \calM^+(D) = \{ y \in \brigplus{,K}[1/t] \otimes_{K_0} D
\mid \phi^n(y) \in \fil^0(K_\infty(\!(t)\!) \otimes_{K_\infty} D_\infty) 
\ \forall n \in \ZZ \}. \]
Si $a,b \in \ZZ$ sont tels que $\fil^{-a+1} D =0$ 
et $\fil^{-b} D = D$ (c'est-à-dire que les \emph{poids} de $D$ sont
dans l'intervalle $[a;b]$), alors on a 
$t^b \cdot \brigplus{,K} \otimes_{K_0} D \subset \calM^+(D) 
\subset t^a \cdot \brigplus{,K} \otimes_{K_0} D$.
\end{rema}

En appliquant la proposition \ref{pgfilhf} à $\dcris(V)$, on retrouve 
le théorème principal de \cite{CH} généralisé dans \cite{BB}.

\begin{coro}\label{crishf}
Si $V$ est une représentation de $G_K$ 
qui devient cristalline sur une extension $K_n$ de $K$,
alors $V$ est de hauteur finie.
\end{coro}

On peut d'ailleurs se demander quand est-ce que $W(D)$ est
cristalline.

\begin{lemm}\label{gamcrdr}
Si $V$ est une représentation de $\Gamma_K$, alors :
\begin{enumerate}
\item $V$ est cristalline si et seulement si 
$V = \oplus_{j \in \ZZ} V^{\Gamma_K = \chi^j}$;
\item $V$ est de Hodge-Tate si et seulement si
elle est potentiellement cristalline.
\end{enumerate}
\end{lemm}

\begin{proof}
Le (1) est l'objet du lemme 3.4.3 de \cite{BP94}, mais nous en 
donnons une nouvelle démonstration. Pour cela, observons
que $\dsen(V) =  K_\infty \otimes_{\Qp} V$ et donc que
$V$ est de Hodge-Tate si et seulement si $\nabla_V$ est
diagonalisable à valeurs propres entières sur 
$K_\infty \otimes_{\Qp}V$; ceci est 
équivalent à demander qu'il existe $n \geq 0$ tel que
$K_n \otimes_{\Qp} V = 
\oplus_{j \in \ZZ} (K_n
\otimes_{\Qp} V)^{\Gamma_{K_n} = \chi^j}$. Comme
$(K_n \otimes_{\Qp} V)^{\Gamma_{K_n} = \chi^j} =
K_n \otimes_{\Qp} V^{\Gamma_{K_n} = \chi^j}$, et qu'une
représentation d'image finie de $\Gamma_K$ est cristalline si
et seulement si elle est triviale, on en déduit le (1) et le (2).
\end{proof}

\begin{coro}\label{wdiscr}
Si $D$ est un $(\phi,\Gamma_K)$-module filtré sur $K_0$, alors la
$B$-paire construite ci-dessus est cristalline si et seulement si 
$D = \oplus_{j \in \ZZ} D^{\Gamma = \chi^j}$.
\end{coro}

\begin{proof}
Comme $W_e = (\bmax \otimes_{K_0} D)^{\phi=1}$, on a 
$\bmax \otimes_{\be} W_e 
= \bmax \otimes_{K_0} D$ ce qui fait que $W$ 
est cristalline si et seulement si $D$ est cristalline en tant
que représentation de $\Gamma_K$. Le corollaire suit alors
du lemme \ref{gamcrdr} ci-dessus.
\end{proof}

\begin{prop}\label{hftopgf}
Si $\dfont$ est un $(\phi,\Gamma_K)$-module de hauteur finie,
alors il existe un $\phi$-module filtré $D$ sur $K_0$ avec 
action de $\Gamma_K$ tel que $\dfont = \brig{}{,K} 
\otimes_{\brigplus{,K}} \calM^+(D)$.
\end{prop}

\begin{proof}
Soit $\dfont^+$ un $(\phi,\Gamma_K)$-module sur 
$\brigplus{,K}$ tel que $\dfont = \brig{}{,K} 
\otimes_{\brigplus{,K}} \dfont^+$ et $D = \dfont^+/X$.
Soit $\nabla = \log(\gamma) / \log_p(\chi(\gamma))$ pour
$\gamma \in \Gamma_K$ proche de $1$.

Commençons par montrer que si $P \in K_0[T]$ est un polynôme
tel que $P(\nabla) : D(j) \to D(j)$ est
bijectif pour tout $j \geq 1$, alors l'application 
$(\dfont^+)^{P(\nabla)=0} \to D^{P(\nabla)=0}$ 
est bijective. 
Pour cela, nous montrons d'abord que l'application :
\[ (K_0[\![t]\!] \otimes_{\brigplus{,K}} \dfont^+)^{P(\nabla)=0} 
\to D^{P(\nabla)=0} \]
est bijective. 
Le fait qu'elle est injective ne pose pas de
problème : si $y \in K_0[\![t]\!] 
\otimes_{\brigplus{,K}} \dfont^+$ est dans son noyau, 
et si $j$ est un  entier $\geq 1$ 
tel que $y \in t^j K_0[\![t]\!] 
\otimes_{\brigplus{,K}} \dfont^+$, alors $P(\nabla)(y)=0$
dans $t^j K_0[\![t]\!] 
\otimes_{\brigplus{,K}} \dfont^+ / t^{j+1} = D(j)$ et donc
$y \in t^{j+1} K_0[\![t]\!] 
\otimes_{\brigplus{,K}} \dfont^+$, 
ce qui fait en itérant que $y=0$. 
Montrons à présent la surjectivité; si $\overline{z} \in D^{P(\nabla)=0}$, 
alors il existe $z_0 \in \dfont^+$ 
qui relève $\overline{z}$ et
tel que $P(\nabla)(z_0) \in X \dfont^+$
et l'hypothèse selon laquelle $P(\nabla) : D(j) \to D(j)$ est bijectif
pour tout $j \geq 0$ nous permet de construire par récurrence $z_j 
\in z_{j-1} + X^j \dfont^+$ tel que $P(\nabla)(z_j) \in 
X^{j+1} \dfont^+$ et donc $z \in (K_0[\![t]\!] 
\otimes_{\brigplus{,K}} \dfont^+)^{P(\nabla)=0}$ 
relevant $\overline{z}$ ce qui fait que notre application est bien bijective.

Montrons à présent que $(K_0[\![t]\!] \otimes_{\brigplus{,K}} 
\dfont^+)^{P(\nabla)=0} = (\dfont^+)^{P(\nabla)=0}$. Pour cela,
remarquons que $(K_0[\![t]\!] \otimes_{\brigplus{,K}} 
\dfont^+)^{P(\nabla)=0}$ est un $K_0$-espace vectoriel 
de dimension finie (puisqu'il s'injecte dans $D$) stable par $\phi$. Si 
l'on choisit une base de $\dfont^+$ et que l'on appelle $Q$ la matrice
de $\phi$ dans cette base, et si l'on choisit une base de 
$(K_0[\![t]\!] \otimes_{\brigplus{,K}} 
\dfont^+)^{P(\nabla)=0}$ et que l'on appelle $P_0$ la matrice
de $\phi$ dans cette base et $Y$ la matrice de passage entre les deux
bases, alors on a $\phi(Y)Q=P_0Y$ et donc $Y=P_0^{-1} \phi(Y) Q$.
Ecrivons $Q = \sum_{k \geq 0} Q_k t^k$ et 
$Y = \sum_{k \geq 0} Y_k t^k$. Si $M$ est une matrice
à coeffcients dans $K_0$, notons $\vp(M)$ le minimum
des valuations de ses coefficients.
Comme $\brigplus{,K} \subset \Qp \otimes_{\Zp} \Zp[\![t/p]\!]$ 
(n'oublions pas que $X=\exp(t)-1$), on voit qu'il existe un entier $h_1$
tel que $\vp(Q_k) \geq - h_1 - k$ et il existe 
par ailleurs un entier $h_2$ 
tel que $\vp(P_0^{-1}) \geq -h_2$; posons $h=h_1+h_2$.
On déduit alors de l'équation $Y=P_0^{-1} \phi(Y) Q$ que
$Y=P_0^{-1} \phi(P_0^{-1}) \cdot \phi^2(Y) \cdot \phi(Q) Q$ et 
(comme $\vp(p^k Q_k) \geq - h_1$) que : 
\begin{align*} 
\vp(Y_k) & \geq - 2 h + \min_{0 \leq \ell \leq k} 
( 2\ell + \vp(Y_\ell) - \ell)  \\
& = - 2 h + \min_{0 \leq \ell \leq k} 
( \ell + \vp(Y_\ell)). 
\end{align*}
On en déduit par récurrence sur $k\geq 2h$
que si $k \geq 2h$, alors $\vp(Y_k) \geq -2h + 
\min_{0 \leq \ell \leq 2h} \vp(Y_\ell)$ et donc,
comme $\OO_{K_0}[\![t]\!] \subset \OO_{K_0}[\![X/p]\!]$, que :
\[ (K_0[\![t]\!] \otimes_{\brigplus{,K}} 
\dfont^+)^{P(\nabla)=0} =
((K_0  \otimes_{\OO_{K_0}} \OO_{K_0}[\![X/p]\!]) 
\otimes_{\brigplus{,K}} 
\dfont^+)^{P(\nabla)=0}. \]
Pour terminer, on utilise de nouveau le fait que 
$((K_0  \otimes_{\OO_{K_0}} \OO_{K_0}[\![X/p]\!]) 
\otimes_{\brigplus{,K}} 
\dfont^+)^{P(\nabla)=0}$ est un $K_0$-espace vectoriel 
de dimension finie stable par $\phi$, et que si $n \geq 0$, alors :
\[\phi(K_0  \otimes_{\OO_{K_0}} \OO_{K_0}[\![X^{p^n}/p]\!]) 
\subset K_0  \otimes_{\OO_{K_0}} \OO_{K_0}[\![X^{p^{n+1}}/p]\!] \]
pour conclure que $(K_0[\![t]\!] \otimes_{\brigplus{,K}} 
\dfont^+)^{P(\nabla)=0} = (\dfont^+)^{P(\nabla)=0}$
puisque $\brigplus{,K} = \cap_{n \geq 0} K_0  \otimes_{\OO_{K_0}} 
\OO_{K_0}[\![X^{p^n}/p]\!]$.

Revenons à notre espace $D=\dfont^+/X$ et soit $Q$ le polynôme 
minimal de $\nabla$. On voit que l'on peut écrire $Q=PR$ où $P$
est non trivial et a la propriété que $P(X+j) \wedge P(X) = 1$ pour
tout $j \geq 1$, ce qui revient à dire que $P(\nabla) : D(j) \to D(j)$
est bijectif. Si $P \neq Q$, alors soit $j$ le plus petit entier tel que
$P(X+j) \wedge P(X) \neq 1$; on vient de montrer que l'application
$(\dfont^+)^{P(\nabla)=0} \to D^{P(\nabla)=0}$ est un isomorphisme,
et si l'on remplace $\dfont^+$ par $\brigplus{,K} \otimes_{K_0} 
(\dfont^+)^{P(\nabla)=0} + X^j \dfont^+$, alors on a toujours
$\dfont = \brig{}{,K} 
\otimes_{\brigplus{,K}} \dfont^+$  mais on peut prendre $P$ de degré
plus grand. En itérant cette opération un nombre fini de fois, on voit donc
que quitte à remplacer $\dfont^+$ par un sous-module, on peut supposer
que l'application $(\dfont^+)^{Q(\nabla)=0} \to D^{Q(\nabla)=0} = D$ 
est un isomorphisme. 

On fait alors un léger abus de notation, et on pose
$D = (\dfont^+)^{Q(\nabla)=0}$; c'est un $\phi$-module sur $K_0$ avec
action de $\Gamma_K$. Remarquons que 
$\brigplus{,K} \otimes_{K_0} D \subset \dfont^+$ et que le déterminant
de l'inclusion est un idéal de $\brigplus{,K}$ stable par $\Gamma_K$ et
par $\phi$ (car $\phi : D \to D$ est bijectif) ce qui fait, par 
le lemme I.3.2
de \cite{Ber3}, que $\brigplus{,K}[1/t] \otimes_{K_0} D 
= \brigplus{,K}[1/t]  \otimes_{\brigplus{,K}} \dfont^+$.
Par le lemme \ref{detpgh}, le déterminant de
$\phi$ sur $\dfont^+$ 
est de la forme $X^{\alpha_0} Q_1^{\alpha_1} \cdots 
Q_m^{\alpha_m}$. Si $n \geq m$, alors l'application déduite du 
Frobenius:
\[ K_{n+1}[\![t]\!] \otimes^{\iota_{n+1}}_{\brigplus{,K}} \dfont^+
\to K_{n+1}[\![t]\!] \otimes^{\iota_n}_{\brigplus{,K}} \dfont^+ \]
est un isomorphisme et la définition : 
\[ \fil^i(D_\infty) = D_\infty \cap
t^i K_\infty[\![t]\!] \otimes^{\iota_n}_{\brigplus{,K}} \dfont^+ \] 
ne dépend donc pas de $n \geq m$. 
Si $y \in \dfont^+$ et $n \geq n(D) = m$, 
alors $\iota_n(y)  \in  \fil^0(K_n (\!(t)\!) 
\otimes_{K_n} D_n)$ et si $n \geq - n(D)$, 
alors $\phi^n(y) = \iota_{n(D)}
(\phi^{n+n(D)}(y))$ ce qui fait que pour tout
$n \in \ZZ$, on a $\phi^n(y) \in 
\fil^0(K_\infty (\!(t)\!) \otimes_{K_\infty} D_\infty)$ et donc
$\dfont^+ \subset \calM^+(D)$. Il reste à constater
que pour tout $n \geq n(D)$, on a :
\[ K_n[\![t]\!] \otimes^{\iota_n}_{\brigplus{,K}} \dfont^+ =
K_n[\![t]\!] \otimes^{\iota_n}_{\brigplus{,K}} \calM^+(D) \]
ce qui fait que $\dfont  = \brig{}{,K} \otimes_{\brigplus{,K}} \dfont^+ = 
\brig{}{,K} \otimes_{\brigplus{,K}}  \calM^+(D)$. 
\end{proof}

Deux $\phi$-modules filtrés sur $K_0$ avec action de $\Gamma_K$
différents peuvent donner le même $(\phi,\Gamma_K)$-module
sur $\brigplus{,K}$. Par exemple, si on a moralement $D_1=K_0$ 
et $D_2=K_0 \cdot t$, alors $\calM^+(D_1) = \calM^+(D_2)
= \brigplus{,K}$. Cet exemple est représentatif comme le 
montre la proposition ci-dessous :

\begin{prop}\label{isomhf}
Si $D_1$ et $D_2$ sont deux  
$\phi$-modules filtrés sur $K_0$ avec action de $\Gamma_K$,
alors $\calM^+(D_1) = \calM^+(D_2)$ si et seulement s'il
existe un isomorphisme $K_0[t,t^{-1}] \otimes_{K_0} D_1
= K_0[t,t^{-1}] \otimes_{K_0} D_2$ compatible à $\phi$
et $\Gamma_K$, et compatible 
à la filtration quand on étend les scalaires à
$K_\infty(\!(t)\!)$.
\end{prop}

\begin{proof}
La construction de $\calM^+(D)$ montre que ce module ne dépend que de
$\brigplus{,K}[1/t] \otimes_{K_0} D$ et de la filtration sur
$K_\infty(\!(t)\!) \otimes_{K_0} D$ ce qui fait que si les conditions
de la proposition sont vérifiées, alors 
$\calM^+(D_1) = \calM^+(D_2)$. On voit réciproquement
que si $\calM^+(D_1) = \calM^+(D_2)$, alors
$\brigplus{,K}[1/t] \otimes_{K_0} D_1 =
\brigplus{,K}[1/t] \otimes_{K_0} D_2$ et cet isomorphisme
est compatible à la filtration quand on étend les scalaires à
$K_\infty(\!(t)\!)$. Si $G_1$ et $G_2$ sont les matrices d'un
élément $\gamma \in \Gamma_K$ qui n'est pas de torsion, et
si $M$ est la matrice de l'isomorphisme entre
$\brigplus{,K}[1/t] \otimes_{K_0} D_1$ et  
$\brigplus{,K}[1/t] \otimes_{K_0} D_2$ alors on a $\gamma(M)G_1
= G_2 M$, et si l'on écrit $M = \sum_{i \gg -\infty} t^i M_i$, alors on voit 
que l'on a $\chi(\gamma)^i M_i G_1 = G_2 M_i$ et donc que si $M_i 
\neq 0$, alors $\chi(\gamma)^i$ est quotient d'une valeur propre de
$G_1$ par une valeur propre de $G_2$, ce qui n'est possible que pour un nombre fini de valeurs de $i$. On en déduit que $M$ est à coefficients dans
$K_0[t,t^{-1}]$ et donc que $K_0[t,t^{-1}] \otimes_{K_0} D_1
= K_0[t,t^{-1}] \otimes_{K_0} D_2$.
\end{proof}

Afin de terminer la démonstration du théorème D, il nous faut montrer 
que si $\brig{}{,K} \otimes_{\brigplus{,K}} \calM^+(D_1) =  
\brig{}{,K} \otimes_{\brigplus{,K}} \calM^+(D_2)$, alors 
$\calM^+(D_1) = \calM^+(D_2)$. Cela suit de la proposition ci-dessous.

\begin{prop}\label{mpisdf}
Si $\dfont = \brig{}{,K} \otimes_{\brigplus{,K}} \calM^+(D)$, et si
$\dfont^+$ est un $(\phi,\Gamma_K)$-module sur $\brigplus{,K}$
de rang fini et contenu dans $\dfont$, alors $\dfont^+ \subset \calM^+(D)$.
\end{prop}

\begin{proof}
Rappelons que $\dfont[1/t] = \brig{}{,K}[1/t] \otimes_{K_0} D$; comme
$\dfont^+$ est de rang fini $r$, il existe $s>0$ et $h\in\ZZ$ tels que
$\dfont^+ \subset t^{-h} \brig{,s}{,K} \otimes_{K_0} D$. Choisissons des bases
de $\dfont^+$ et de $D$ et appelons $M$, $P^+$ et $P_0$ la matrice de
passage, la matrice de $\phi$ sur $\dfont^+$ et la matrice de $\phi$ sur $t^{-h} D$
ce qui fait que $P_0 \phi(M) = M P^+$. Comme $P_0 \in \mathrm{GL}_d(K_0)$, le corollaire I.4.2 de \cite{Ber3} (qui s'étend verbatim aux matrices rectangulaires) 
nous donne que 
$M \in \mathrm{M}_{d \times r}(\brigplus{,K})$ et donc que 
$\dfont^+ \subset t^{-h} \brigplus{,K} \otimes_{K_0} D$.

Quitte à augmenter $s$, on a de plus que 
$\dfont^+ \subset \brig{,s}{,K} \otimes_{\brigplus{,K}} \calM^+(D)$ 
et si $y \in \dfont^+$ et $n \geq n(s)$, alors $\iota_n(y) \in 
\fil^0(K_\infty (\!(t)\!) \otimes_{K_\infty} D_\infty)$ 
puisque $\iota_n(\brig{,s}{,K}) \subset K_n[\![t]\!]$,
et donc $\iota_n(\dfont^+) \subset 
\fil^0(K_\infty (\!(t)\!) \otimes_{K_\infty} D_\infty)$. Enfin, si 
$n \geq  - n(s)$ et $y \in \dfont^+$, alors
$\phi^n(y) = \iota_{n(s)}(\phi^{n(s)+n}(y))$ et donc 
$\phi^n(\dfont^+) \subset 
\fil^0(K_\infty (\!(t)\!) \otimes_{K_\infty} D_\infty)$ ce qui fait,
étant donnée la définition de $\calM^+(D)$ donnée dans la 
proposition \ref{pgfilhf}, que $\dfont^+ \subset \calM^+(D)$.
\end{proof}

En particulier, un $(\phi,\Gamma_K)$-module sur $\brigplus{,K}$
de rang fini et contenu dans $\dfont$ est nécessairement de rang $\leq d$.

\appendix

\section{Liste des notations}

Voici une liste des principales notations dans l'ordre où elles apparaissent. 

Introduction : $k$, $K$, $K_0$, $G_K$.

\S 1.1 : $\et^+$, $\eps^{(n)}$,
$\et$, $\atplus$, $\atfont$, $\btplus$, $\bt$, $\theta$,
$\bdr^+$, $X$, $t$, $\bdr$, $\tilde{p}$, $\bmax^+$, $\btp$,
$\bmax$, $\bst$, $N$, $\bcris$, $\btdag{,r}$, 
$\btrig{,r}$, $K_n$, $H_K$, $\Gamma_K$, $\brig{,r}{,K}$,
$K_0'$, $\bdag{,r}_K$, $M_{a,h}$, $\be$, $t_h$, $V_{[r;s]}$,
$\bt^{[r;s]}$, $\bfont^{[r;s]}_K$, $\rho(r)$.

\S 1.2 : $\deg(\dfont)$, $\NP(\dfont)$, $r_n$, $n(r)$, $\iota_n$, $Q_n$.   

\S 2 : $W_e$, $W_{dR}^+$, $W_{dR}$.

\S 2.2 : $\dtilde^r(W)$, $W(\dfont)$.

\S 2.3 : $\dcris(W)$, $\dst(W)$, $\ddr(W)$, $\calM(D)$, $\dsen(W)$.

\S 3.2 : $\rep(a,h)$.

\S 3.3 : $\brigplus{,K}$, $\bplus_K$, $\calM^+(D)$, $\nabla$.

\end{document}